\documentclass[reqno]{amsart}

\usepackage{fancyhdr} 
\usepackage[lite]{amsrefs} 
\usepackage{amssymb} 
\usepackage[all,cmtip]{xy}
\usepackage{graphicx} 
\usepackage{euscript}
\usepackage{color}
\usepackage{array}
\usepackage{picture}
\usepackage{epic}
\usepackage{tikz}
 \usetikzlibrary{backgrounds,cd}
\usepackage{comment}
\usepackage{enumerate}
\usepackage{dynkin-diagrams}
\usepackage{hyperref}

\newcommand{\isom}{\cong}
\newcommand{\set}[1]{\left\{#1\right\}} 
\newcolumntype{H}{>{\setbox0=\hbox\bgroup}c<{\egroup}@{}} 

\DeclareMathOperator{\codim}{codim}
\DeclareMathOperator{\image}{Im}
\DeclareMathOperator{\pr}{pr}
\DeclareMathOperator{\id}{id}

\DeclareMathOperator{\Gr}{Gr}

\newcommand{\cO}{\mathcal O}

\newcommand{\bC}{\mathrm{\bf C}}           

\newcommand{\bQ}{\mathrm{\bf Q}}           

\newcommand{\bP}{\mathrm{\bf P}}          

\theoremstyle{plain} 
\newtheorem{theorem}[equation]{Theorem}
\newtheorem{corollary}[equation]{Corollary}
\newtheorem{lemma}[equation]{Lemma}
\newtheorem{key lemma}[equation]{Key Lemma}
\newtheorem{proposition}[equation]{Proposition}
\newtheorem{facts}[equation]{Facts}

\theoremstyle{definition}
\newtheorem{definition}[equation]{Definition}

\theoremstyle{remark}
\newtheorem{remark}[equation]{Remark}
\newtheorem{example}[equation]{Example}

\begin{document}


\title{Decompositions of Schubert Varieties\\ and Small Resolutions}

\author{Scott Larson}





\date{\today}
\keywords{small resolutions, Schubert varieties, flag varieties, Bott-Samelson, BP decomposition, fiber bundles, pattern avoidance}
\subjclass[2010]{22E47, 14M15}

\begin{abstract}
We provide a method for gluing (small) resolutions of singularities of Schubert varieties \(X_w\).
An explicit isomorphism of \(X_w\) with an (iterated) bundle is constructed when \(w\) has an (iterated) BP decomposition.
Combined with the first result this gives many new small resolutions of Schubert varieties.
In type A, this can be expressed in terms of pattern avoidance.
Also we show resolutions of Schubert varieties constructed quite generally are in fact Gelfand-MacPherson resolutions.
\end{abstract}

\maketitle
\thispagestyle{empty}

\setcounter{tocdepth}{1}
\tableofcontents
\section*{Introduction}

Geometry of Schubert varieties governs certain key properties of the representation theory in category \(\cO\) and of complex reductive groups.
The importance of Schubert varieties is reflected by vast literature on their geometry, including a rich history of their desingularizations.
Small resolutions of singularities can enable computing Kazhdan-Lusztig polynomials (as in, e.g., \cite{Zelevinskii83}) and characteristic cycles (as in \cite{BresslerFinkelbergLunts90}).

Our original motivation for this paper was to construct small resolutions.
Let \(G\) be a complex connected reductive algebraic group and \(B\) a fixed Borel subgroup.
For notational convenience, we let \(G_w=\overline{B\dot wB}\subseteq G\), where quotienting by \(B\) gives a Schubert variety in the flag variety \(G/B\).
It was soon realized that all resolutions of Schubert varieties we studied from the literature are particular examples of the morphism 
\begin{equation}\label{equation: mu in G}
\mu\,\colon G_{w_0}\times^{R_1}\cdots\times^{R_m}G_{w_m}\to G_w
\end{equation}
(as described in Definition~\ref{definition: mu}), where \(R_i\) is a parabolic subgroup stabilizing \(G_{w_{i-1}}\) and \(G_{w_i}\), and \(\mu\) is defined by multiplication.
This suggests that if any \(G_w/B\) has a small resolution, then \(G_w/B\) has a small resolution of the form \eqref{equation: mu in G} (after quotienting by \(B\) on the right).
Moreover, in our quest for finding small fiber dimensions, it was eventually realized that all fiber bundle structures on Schubert varieties can be described by the same formula -- namely, when fibers of \(\mu\) have dimension zero.

The morphism \(\mu\) will always be \((B\times B)\)-equivariant, but it is a recurring obstacle to check for equivariance with respect to the stabilizer of \(G_w\) in \(G\times G\).
When the map
\begin{equation}\label{equation: m is one}
\mu\,\colon G_v\times^RG_w\to G_u
\end{equation}
is an isomorphism (as described in Corollary~\ref{corollary: fiber bundle decomp}), it is interesting that the natural action on \(G_u\) is typically larger than the action on \(G_v\times^RG_w\).

Let \(R\subseteq P\) be parabolic subgroups of \(G\) containing \(B\).
Richmond-Slofstra \cite{RichmondSlofstra16} describe exactly when the morphism \(\pi\,\colon G/R\to G/P\) restricts to a fiber bundle on \(G_wR/R\), by a Coxeter-theoretic condition called BP decomposition.
Our main result in this direction describes the fiber bundle structure on \(G_wR/R\) explicitly as a Bott-Samelson type variety by using \eqref{equation: m is one}, as in Proposition~\ref{proposition: BP isom}.

Our key Lemma~\ref{lemma: small glue} shows how to take two small resolutions of the form \eqref{equation: mu in G} and construct new resolutions using \eqref{equation: m is one}.
Applying the lemma to isomorphisms of the form \eqref{equation: m is one} requires having enough equivariance, and what we are often able to show is that a small resolution \eqref{equation: mu in G} of \(G_w\) can be made maximally equivariant, satisfying
\begin{subequations}\label{equation: equivariance}
\begin{align}
\tau(w^{-1})&=\tau(w_0^{-1}) \label{equation: a}\\
\tau(w)&=\tau(w_m) \label{equation: b},
\end{align}
\end{subequations}
where \(\tau\) is defined in \eqref{equation: tau}.
These conditions guarantee that any standard parabolic subgroup of \(G\) stabilizing \(G_w\) by left multiplication also stabilizes \(G_{w_0}\) on the left, and likewise for \(G_w\) and \(G_{w_m}\) on the right.
Equations~\eqref{equation: equivariance} are indeed stronger than necessary for applying Lemma~\ref{lemma: small glue}.

Using the above methods we obtain our goal of explicitly constructing new small resolutions for families of Schubert varieties, e.g., in Proposition~\ref{proposition: small BP family}.
This family is best described using a pattern avoidance result of \cite{AllandRichmond18}.
We view Lemma~\ref{lemma: small glue} as highlighting the importance of determining small resolutions in low rank.
On the way, we classify all Schubert varieties for \(W\) of type \(A_{n-1}\) (\(n\leq6\)) admitting \emph{any} small resolution.

It is important to describe small resolutions explicitly, and in some cases multiple nonisomorphic small resolutions may occur (as is well-known).
In particular, fiber dimensions are needed to determine whether the resolution is small, and then cohomology of fibers are needed to compute intersection cohomology.
Thus in our opinion, another main result of this paper appears in \S\ref{section: Gelfand-MacPherson}, which simplifies a large family of resolutions to the form constructed by Gelfand-MacPherson \cite{GelfandMacPherson82} for which there exists a formula for \(\dim(\mu^{-1}(\mathrm{pt}))\).
This relies essentially on \cite{RichmondSlofstra16}, and our perspective on the corresponding results.
In particular, every smooth Schubert variety \(G_w\) in a simply laced group admits an isomorphism of the form \eqref{equation: mu in G}, where \eqref{equation: equivariance} holds true, such that for every \(0\le i\le m\), \(G_{w_i}\) is a parabolic subgroup. 

I thank Edward Richmond for explaining the remarkable work in Richmond-Slofstra \cite{RichmondSlofstra16} to me, and I thank Roger Zierau for many helpful conversations.
\numberwithin{equation}{section}

\section{Bott-Samelson Type Varieties}

Bott-Samelson \cite{BottSamelson58} constructed certain quotients by group actions in the category of smooth manifolds, which proved useful in studying the topology of compact Lie groups and symmetric spaces.
The same construction for algebraic varieties has been useful in studying properties of Schubert varieties.
We recall this construction here and apply it to Schubert varieties.
Then we describe a proper map \(\mu\) from such a variety to a flag variety.
We conclude this section by characterizing when \(\mu\) is a resolution of singularities of a Schubert variety.
Of particular interest is characterizing when \(\mu\) is birational and when \(\mu\) is an isomorphism.

Let \(X\) be an algebraic variety and let \(H\) be a linear algebraic group.
Suppose that \(X\) is a \(H\)-variety with a right action.
Let \(X/H\) be the quotient space with the quotient topology, let \(\pi\,\colon X\to X/H\) be the quotient map, and for any \(U\subseteq X/H\) open let \(\cO_{X/H}(U)\) be the set of functions \(f\,\colon U\to k\) such that \(f\circ\pi\vert\pi^{-1}(U)\) is in \(\cO_X(\pi^{-1}(U))\).
Thus, \(\cO_{X/H}(U)\) may be identified with the ring of invariant functions \(\cO_X(\pi^{-1}(U))^H\) on \(\pi^{-1}(U)\).
Then \(X/H\) is a ringed space, but may fail to be an algebraic variety.
All quotients we consider will be varieties, in particular they occur naturally as subvarieties of a quotient \(G\times^{H_1}\cdots\times^{H_m}G\) (as in \eqref{equation: Bott-Samelson variety}).

Suppose that the right action of \(H\) on \(X\) is free and let \(Y\) be a left \(H\)-variety.
Then \(X\times Y\) is a \(H\)-variety with a free right action by \((x,y)h=(xh,h^{-1}y)\).
Let \(X\times^HY\) denote the quotient space \((X\times Y)/H\) and let \(\rho\,\colon X\times Y\to X\times^HY\) be the quotient map.
There exists a natural map of \(X\times^HY\) onto \(X/H\) which makes the following diagram commutative:
\[\begin{tikzcd}
        X\times Y\arrow{r}{}\arrow{d}{}&X\arrow{d}{}\\
        X\times^HY\arrow{r}{}&X/H
\end{tikzcd}\]
where the other maps are the natural quotient maps.

One may check the quotient \(X\times^HY\times^{H'}Y'\) is isomorphic to both \((X\times^HY)\times^{H'}Y'\) and \(X\times^H(Y\times^{H'}Y')\) as ringed spaces.
As noted in \cite{BottSamelson58}, there is an obvious extension to more factors to obtain ringed spaces 
\begin{equation}\label{equation: Bott-Samelson quotients}
X_0\times^{H_1}X_1\times^{H_2}\cdots\times^{H_m}X_m.
\end{equation}

Let \(G\) be an algebraic group and let \(H_1,\ldots,H_m\) be closed subgroups.
For every \(1\le i\le m\), \(H_i\) acts freely on \(G\) by multiplication on the right and \(G/H_i\) is an algebraic variety.
Observe that for every \(1\le i<m\), \(H_i\) also acts on \(G\) by multiplication on the left, such that the actions of \(H_i\) and \(H_{i+1}\) on \(G\) associate.
The quotient \(G\times^{H_1}\cdots\times^{H_m}G\) is an algebraic variety.

To see this, define \(\tilde\varphi\,\colon G\times\cdots\times G\to G\times\cdots\times G\) by
\begin{equation}\label{equation: phi in G}
\tilde\varphi(g_0,\ldots,g_m)=(g_0,g_0g_1,\ldots,g_0g_1\cdots g_m)
\end{equation}
a morphism of varieties.
Let \(\pi\,\colon G\times\cdots\times G\to G/H_1\times\cdots\times G/H_m\times G\) denote the projection morphism of varieties.
Then \(\pi\circ\tilde\varphi\) is constant on \((H_1\times\cdots\times H_m)\)-orbits, so gives a morphism of ringed spaces \(\varphi\,\colon G\times^{H_1}\cdots\times^{H_m}G\to G/H_1\times\cdots\times G/H_m\times G\), where
\begin{equation}\label{equation: phi}
\varphi[g_0,\ldots,g_m]=(g_0H_1/H_1,\ldots,g_0\cdots g_{m-1}H_m/H_m,g_0\cdots g_m).
\end{equation}

Define \(\tilde\psi\,\colon G\times\cdots\times G\to G\times\cdots\times G\) by
\begin{equation}\label{equation: psi in G}
\tilde\psi(g_0,\ldots,g_m)=(g_0,g_0^{-1}g_1,g_1^{-1}g_2,\ldots,g_{m-1}^{-1}g_m),
\end{equation}
a morphism of varieties.
Note \(\tilde\psi\) is the inverse morphism of \(\tilde\varphi\).
Let \(\rho\,\colon G\times\cdots\times G\to G\times^{H_1}\cdots\times^{H_m}G\) be the quotient morphism.
Then \(\rho\circ\tilde\psi\) is constant on \((H_1\times\cdots\times H_m)\)-orbits, so gives a morphism of ringed spaces \(\psi\,\colon G/H_1\times\cdots\times G/H_m\times G\to G\times^{H_1}\cdots\times^{H_m}G\), where
\begin{equation*}\label{equation: psi}
\psi(g_1H_1/H_1,\ldots,g_{m-1}H_m/H_m,g_m)=[g_1,g_1^{-1}g_2,g_2^{-1}g_3,\ldots,g_{m-1}^{-1}g_m].
\end{equation*}
Then \(\varphi\) and \(\psi\) are inverse morphisms, so the ringed spaces are isomorphic.
Hence the quotient is an algebraic variety.

All of our quotients will embed naturally into \(G\times^{H_1}\cdots\times^{H_m}G\) and again be algebraic varieties, so we can think of this as providing a safety zone for the ringed spaces on the quotient to be an algebraic variety.
We illustrate this point of view by considering a simple case before presenting the general construction in \eqref{equation: Bott-Samelson variety}.

Let \(G_0,\ldots,G_m\) be closed subgroups of \(G\), and for every \(1\le i\le m\), let \(H_i\) be a closed subgroup of \(G\) in \(G_{i-1}\cap G_i\).
Then \(G_0\times^{H_1}\cdots\times^{H_m}G_m\) is an algebraic variety.
Indeed, the morphism \(\tilde\varphi\) defined in \eqref{equation: phi in G} restricts to a closed embedding \(\tilde\iota\,\colon G_0\times\cdots\times G_m\to G\times\cdots\times G\).
Hence the morphism \(\tilde\psi\) defined in \eqref{equation: psi in G} restricts to the inverse morphism, so \(\iota\) is an embedding of ringed spaces.
The image of \(\iota\) is the fibered product
\begin{equation}\label{equation: iterated fibered product}
G_0/H_1\underset{G/G_1}{\times}G/H_2\underset{G/G_2}{\times}\cdots\underset{G/G_{m-1}}{\times}G/H_m\underset{G/G_m}{\times}G
\end{equation}
which is a closed subvariety of the product.
Therefore \(\iota\) is a closed embedding of algebraic varieties.


Our main construction uses Schubert varieties to construct quotient varieties as in \eqref{equation: Bott-Samelson quotients}.
We note that this idea has appeared before by various authors (as in, e.g., \cite{RichardsonSpringer90}).
The resulting varieties will be iterated fiber bundles of the corresponding Schubert varieties.

From now on, let \(G\) be a connected reductive algebraic group and fix a Borel subgroup \(B\) along with a maximal torus \(T\subseteq B\).
Let \(X=G/B\) be the flag variety of \(G\).
There are finitely many \(B\)-orbits on \(X\)
\begin{equation*}\label{equation: Bruhat decomposition}
X=\coprod_{w\in W}B\dot wB/B
\end{equation*}
where \(W=N_G(T)/T\) is the Weyl group of \(G\) and \(\dot w\) is a representative of \(w\) in \(N_G(T)\), the normalizer of \(T\).
Let \(S\) be the set of simple reflections in \(W\) with respect to \(B\).

We identify a standard parabolic subgroup \(B\subseteq P_I\subseteq G\) with subsets of simple reflections \(\emptyset\subseteq I\subseteq S\) such that for every \(s\in I\), \(P_I\) contains \(\dot s\) in \(G\).
In particular, \(P_I\) has semisimple rank \(\#I\).
We write \(X^I=G/P_I\) for the flag variety of \(G\) corresponding to \(I\).
There are finitely many \(B\)-orbits on \(X^I\)
\begin{equation*}\label{equation: generalized Bruhat decomposition}
X^I=\coprod_{w\in W^I}B\dot wP_I/P_I
\end{equation*}
where \((W_I,I)\) is the Weyl group generated by \(I\), and we let \(W^I\) be the set of maximal length representatives of cosets in \(W/W_I\).
We set \(w_I=\max(W_I)\), so in particular, \(w_S=\max(W)\).

Given \(w\in W\), let
\begin{equation*}\label{equation: Schubert varieties}
\begin{split}
G_w&=\overline{B\dot wB}\\
X_w&=\overline{B\dot wB/B}\\
X_w^I&=\overline{B\dot wP_I/P_I}
\end{split}
\end{equation*}
be closures in \(G\), \(X\), and \(X^I\).
In particular, \(G_{w_I}=P_I\).
Let \(\pi\,\colon G\to X\) be the quotient map, which is a fiber bundle with fiber \(B\).
Then base change of \(\pi\) with respect to the inclusion \(X_w\subseteq X\) gives a fiber bundle \(\pi'\,\colon G_w\to X_w\) with fiber \(B\) (as shown in Lemma~\ref{lemma: base change}).

The main construction used in this paper is as follows.
Let \(w_0,\ldots,w_m\in W\).
If for every \(1\le i\le m\), \(H_i\) stabilizes \(G_{w_{i-1}}\) by right multiplication and \(G_{w_i}\) by left multiplication, the quotient
\begin{equation}\label{equation: Bott-Samelson variety}
G_{w_0}\times^{H_1}\cdots\times^{H_m}G_{w_m}
\end{equation}
is a well-defined ringed space -- and in fact is an algebraic variety.
Indeed, define \(\iota\,\colon G_{w_0}\times^{H_1}\cdots\times^{H_m}G_{w_m}\to G/H_1\times\cdots\times G/H_m\times G\) by
\begin{equation}\label{equation: iota}
\iota[g_0,\ldots,g_m]=(g_0H_1/H_1,g_0g_1H_2/H_2,\ldots,g_0\cdots g_m)
\end{equation}
the same formula as in \eqref{equation: phi}.
The diagram given by universal properties of quotients
\[\begin{tikzcd}
        G_{w_0}\times\cdots\times G_{w_m}\arrow{r}{\subseteq}\arrow{d}{\rho'}&G\times\cdots\times G\arrow{d}{\rho}\\
        G_{w_0}\times^{H_1}\cdots\times^{H_m}G_{w_m}\arrow{r}{}&G\times^{H_1}\cdots\times^{H_m}G
\end{tikzcd}\]
commutes.
It follows that the image \(Z\) of \(\iota\) is closed, since the closed subset \(G_{w_0}\times\cdots\times G_{w_m}\) of \(G\times\cdots\times G\) is \(\rho\)-saturated and \(\rho\) is a surjective open map of topological spaces.
Similar to \eqref{equation: iterated fibered product}, \(\iota\) is a closed embedding of algebraic varieties.

An algebraic variety of the form \eqref{equation: Bott-Samelson variety} was considered independently by \cite{Demazure74} and \cite{Hansen73}, with \(G_{w_i}\) minimal parabolic subgroups (\(w_i=s_i\in S\)), and \(H_i\) Borel subgroups.
This variety enjoys many nice properties, such as being a smooth iterated fiber bundle of \(\bP^1\)'s, after quotienting by a Borel subgroup on \(G_{w_m}\).
This construction was used to provide a resolution of singularities for any Schubert variety.
Resolutions of this form are often called Demazure resolution, Bott-Samelson resolution, Bott-Samelson-Demazure-Hansen resolution, etc.

Demazure's resolution was generalized by Gelfand-MacPherson \cite{GelfandMacPherson82} using more general parabolic subgroups, i.e., \(w_i=w_{I_i}\), for various \(I_i\subseteq S\).
The quotients we consider \eqref{equation: Bott-Samelson variety} may be viewed as generalizations of those of \cite{GelfandMacPherson82}.
These will have an iterated fiber bundle structure, but will not in general be smooth.

If for every \(1\le i\le m\), we take \(H_i\) to be a parabolic subgroup \(R_i\), then we define a proper algebraic morphism \(\mu\,\colon G_{w_0}\times^{R_1}\cdots\times^{R_m}G_{w_m}/B\to X\) and show the domain is an iterated fiber bundle.
We will be interested in cases where \(\mu\) is either an isomorphism or is a resolution of singularities (onto its image), so we provide a proof that \(\mu\) is always proper, describe precisely when \(\mu\) is birational, and characterize when the domain of \(\mu\) is smooth in terms of Weyl group elements.

\begin{definition}\label{definition: mu}
Define \(\mu\,\colon G_{w_0}\times^{R_1}\cdots\times^{R_m}G_{w_m}/B\to X\) by
\begin{equation}\label{equation: mu}
\mu[g_0,\ldots,g_mB/B]=g_0\cdots g_mB/B.
\end{equation}
\end{definition}

Observe that for every \(0\le i\le m\), \(B\) stabilizes \(G_{w_i}\) on the left and on the right.
From here and below, we assume that for every \(1\le i\le m\), the parabolic subgroup \(R_i\) is a parabolic subgroup containing \(B\) corresponding to simple reflections \(J_i\).

\begin{proposition}\label{proposition: mu}
\

\begin{enumerate}[(i)]
\item The map \(\mu\,\colon G_{w_0}\times^{R_1}\cdots\times^{R_m}G_{w_m}/B\to X\) defined in \eqref{equation: mu} is a proper algebraic morphism with image \(X_w\) for some \(w\in W\).
\item The variety \(G_{w_0}\times^{R_1}\cdots\times^{R_m}G_{w_m}/B\) is an iterated Zariski locally trivial fiber bundle.
\item The map \(\mu\,\colon G_{w_0}\times^{R_1}\cdots\times^{R_m}G_{w_m}/B\to X_w\) is birational if and only if 
\begin{equation}\label{equation: decomposition}
\ell(w)=\sum_{i=0}^m\ell(w_i)-\sum_{i=1}^m\ell(w_{J_i})
\end{equation}
i.e., \(\mu\) is birational if and only if it is generically finite.
\item The variety \(G_{w_0}\times^{R_1}\cdots\times^{R_m}G_{w_m}/B\) is smooth if and only if for every \(0\le i\le m\), the Schubert variety \(X_{w_i}\) is smooth.
\end{enumerate}

\begin{proof}[Proof of (i)]
Define \(\iota\,\colon G_{w_0}\times^{R_1}\cdots\times^{R_m}G_{w_m}/B\to X^{J_1}\times\cdots\times X^{J_m}\times X\) by
\begin{equation*}\label{equation: iota mod B}
\iota[g_0,\ldots,g_mB/B]=(g_0R_1/R_1,g_0g_1R_2/R_2,\ldots,g_0\cdots g_mB/B).
\end{equation*}
Similar to \eqref{equation: iota}, \(\iota\) is a closed embedding.
Hence \(G_{w_0}\times^{R_1}\cdots\times^{R_m}G_{w_m}/B\) is a projective variety.
It follows that \(\mu\) is a proper algebraic morphism, since any algebraic morphism between projective varieties is a projective morphism.
The quotient \(G_{w_0}\times^{R_1}\cdots\times^{R_m}G_{w_m}/B\) is irreducible since it is the image of the irreducible product \(G_{w_0}\times\cdots\times{G_{w_m}}\) under the quotient morphism.
It follows that the image of \(\mu\) is closed (\(\mu\) is proper), irreducible (the domain of \(\mu\) is irreducible), and \(B\)-stable (\(\mu\) is \(B\)-equivariant).
Hence the image of \(\mu\) is equal to \(X_w\) for some \(w\in W\).
\end{proof}

\begin{proof}[Proof of (ii)]
Let \(R\) be a standard parabolic subgroup corresponding to \(J\subseteq S\), and let \(w\in W\) such that \(R\) stablizes \(G_w\) by right multiplication.
By \cite{Jantzen03}, the map \(\pi\,\colon G\to X^J\) has local sections.
It follows that the base change to \(G_w\to X_w^J\) also has local sections.
By \cite[\S 5.5.8]{Springer09}, the map \(G_{w_0}\times^{R_1}Y\to X_{w_0}^{J_1}\) has local sections, where \(Y=G_{w_1}\times^{R_2}\cdots\times^{R_m}G_{w_m}/B\).
The claim follows by recursion.
\end{proof}

\begin{proof}[Proof of (iii)]
By part (ii), the dimension of the domain of \(\mu\) is
\begin{equation*}\label{equation: mu dimension}
\dim(X_{w_0}^{J_1})+\dim(Y)=\ell(w_0)-\ell(w_{J_1})+\dim(Y)
\end{equation*}
where the equality follows from the fiber bundle \(X_{w_0}\to X_{w_0}^{J_1}\) with fiber \(R_1/B\).
Hence the dimension of the domain of \(\mu\) is \(\sum_{i=0}^m\ell(w_i)-\sum_{i=1}^m\ell(w_{J_i})\) by recursion.
It follows that \(\mu\) is generically finite if and only if \eqref{equation: decomposition} holds true.

It remains to show that \(\mu\) is in fact birational.
We will first show that
\begin{equation*}\label{equation: generic fiber direct product}
\mu^{-1}(B\dot wB/B)\isom B\dot wB/B\times \mu^{-1}(\dot wB/B).
\end{equation*}
Then \eqref{equation: decomposition} will force \(\mu^{-1}(\dot wB/B)\) to be finite by generic finiteness.
The irreducibility of the domain of \(\mu\) forces the open set \(\mu^{-1}(B\dot wB/B)\) to be irreducible, which will then force \(\mu^{-1}(\dot wB/B)\) to be a single point.

Let \(U\) be the unipotent radical of \(B\).
The closed irreducible subgroup \(U_{w^{-1}}=U\cap\dot w\dot w_SU\dot w_S\dot w^{-1}\) of \(U\) is described by \cite[\S 8.3.5]{Springer09}.
Then \cite[\S 8.3.6]{Springer09} gives an isomorphism \(\eta\,\colon B\dot wB/B\isom U_{w^{-1}}\) having inverse \(u\mapsto u\dot wB/B\).
Define \(\alpha\,\colon B\dot wB/B\times\mu^{-1}(x)\to\mu^{-1}(B\dot wB/B)\) by \(\alpha(x,y)=\eta(x)y\), where \(x\in B\dot wB/B\) and \(y\in\mu^{-1}(\dot wB/B)\).

Define \(\beta\,\colon\mu^{-1}(B\dot wB/B)\to B\dot wB/B\times\mu^{-1}(\dot wB/B)\) by \(\beta(y)=(\mu(y),\eta(\mu(y))^{-1}y)\), where \(y\in\mu^{-1}(B\dot wB/B)\).
It is a routine calculation to check that \(\alpha\) and \(\beta\) are inverse regular maps.
It follows that \(\mu\) is birational.
\end{proof}

\begin{proof}[Proof of (iv)]
By part (ii), the variety \(G_{w_0}\times^{R_1}\cdots\times^{R_m}G_{w_m}/B\) is the total space of a fiber bundle with base \(X_{w_0}^{J_1}\) and fiber \(Y=G_{w_1}\times^{R_2}\cdots\times^{R_m}G_{w_m}/B\).
Hence the fiber bundle is smooth if and only if both \(X_{w_0}^{J_1}\) and \(Y\) are smooth.

The fiber bundle \(X_{w_0}\to X_{w_0}^{J_1}\) has smooth fiber \(R_1/B\) since it is the base change of \(X\to X^{J_1}\) by \(X_{w_0}^{J_1}\subseteq X^{J_1}\) inclusion.
Hence \(X_{w_0}^{J_1}\) is smooth if and only if \(X_{w_0}\) is smooth.
By recursion, the iterated fiber bundle is smooth if and only if for every \(0\le i\le m\), the Schubert variety \(X_{w_i}\) is smooth.
\end{proof}
\end{proposition}

\section{The Monoid $(W,\star)$}
We use \(\mu\) to define a monoid product \(\star\) on the Weyl group \(W\).
The monoid is used to define a function \(\tau\) from \(W\) to subsets of simple reflections, where \(\tau(w)\) is often called the \(\tau\)-invariant or the right descents of \(w\).
Comparing \(\tau\) to another function \(\sigma\) called the support will be used repeatedly in this paper to study \((W,\star)\).
This monoid coincides with that of Richardson-Springer \cite{RichardsonSpringer90}.


Recall that a \emph{monoid} is a set together with an associative law of composition \(M\times M\to M\) such that \(M\) contains an identity element.

Define \(\star\,\colon W\times W\to W\) by \(v\star w=u\), where \(X_u\) is the image of \(\mu\,\colon G_v\times^RX_w\to X\) as in Proposition~\ref{proposition: mu} part (i).
Explicitly, we have
\begin{equation*}\label{equation: mu image}
X_{v\star w}=\image(\mu)=\overline{B\dot vB\dot wB/B}
\end{equation*}
where \(\mu\) is defined by any parabolic subgroup \(R\supseteq B\) such that \(R\) stabilizes \(G_v\) on the right and \(X_w\) on the left (the image of \(\mu\) does not depend on the choice of \(R\)).
Equivalently, \(v\star w\) may be defined in \(G\) by
\begin{equation}\label{equation: mu image in G}
G_{v\star w}=G_vG_w
\end{equation}
i.e., \(\overline{B(\dot{v\star w})B}=\overline{B\dot vB\dot wB}\).
The monoid associativity and identity element properties are easy to see.

\begin{facts}\label{facts: monoid}
There are a few easy facts about \(\star\) that we will use many times.
\begin{enumerate}[(a)]
\item For any \(v,w\in W\), we have \((v\star w)^{-1}=w^{-1}\star v^{-1}\).
\item If \(J\subseteq S\) and \(w_J\) is the long element of \(W_J\), then \(G_{w_J}=P_J\) and \(G_{v\star w_J}=G_vP_J\).
\item \(v,w,vw\leq v\star w\).
\end{enumerate}

\begin{proof}
(a) Apply the group inverse to \eqref{equation: mu image in G}.
(b) See \cite[\S 8.4.3]{Springer09}.
(c) We have \(X_v,X_w,X_{vw}\subseteq G_vG_w/B=X_{v\star w}\) since \(B\) is contained in each \(G_u\).
\end{proof}
\end{facts}


Consider \(w\in W\) and \(s\in S\).
It is well-known that
\begin{equation*}\label{equation: cell multiplication}
G_wG_s=\overline{B\dot wB\dot sB}=
\begin{cases}
        \overline{B\dot w\dot sB},&\ell(ws)=\ell(w)+1,\\
        \overline{B\dot wB},&\ell(ws)=\ell(w)-1.
\end{cases}
\end{equation*}
e.g., by \cite[\S8.3.7]{Springer09}.
Therefore,
\begin{equation}\label{equation: monoid simples}
w\star s=
\begin{cases}
        ws,&\ell(ws)=\ell(w)+1,\\
        w,&\ell(ws)=\ell(w)-1.
\end{cases}
\end{equation}

\begin{definition}
Define the function \(\tau\) from \(W\) to subsets of \(S\) by
\begin{equation*}\label{equation: tau}
\tau(w)=\set{s\in S\mid w\star s=w},
\end{equation*}
called the \emph{\(\tau\)-invariant} of \(w\).
This set of simple reflections is often called the \emph{right descents} of \(w\).
The terminology of \(\tau\)-invariant is popular in representation theory (as, e.g., in \cite{BorhoBrylinski85}).
\end{definition}

\begin{facts}\label{facts: tau}
The following is a list of easy facts.
\begin{enumerate}[(a)]
\item \(\tau(w^{-1})=\set{s\in S\mid s\star w=w}\).
\item If \(s_\alpha\) is the reflection in the simple root \(\alpha\), then \(s_\alpha\) is in \(\tau(w)\) if and only if \(w\alpha<0\).
\item \(\tau(w)\subseteq\tau(v\star w)\).
\item If \(J\subseteq S\) then \(\tau(w_J)=J\).
\end{enumerate}

\begin{proof}
(a) The relation \(\overline{B\dot sB\dot wB}^{-1}=\overline{B\dot w^{-1}B\dot sB}\) holds in \(G\).
(b) Let \(\Phi\) be the set of roots of \((G,T)\), let \(\Phi^+\) be the roots of \((B,T)\), and let
\begin{equation*}\label{equation: Phi(w)}
\Phi(w)=\set{\alpha\in\Phi^+\mid w\alpha\in-\Phi^+}
\end{equation*}
be the \emph{right inversions} of \(w\).
Then
\begin{equation*}\label{equation: Phi(ws)}
\begin{split}
\Phi(ws_\alpha)=
\begin{cases}
s_\alpha\Phi(w)\cup\set{\alpha},&\text{ if }w\alpha\in\Phi^+,\\
s_\alpha(\Phi(w)\smallsetminus\set{\alpha}),&\text{ if }w\alpha\in-\Phi^+,
\end{cases}
\end{split}
\end{equation*}
by \cite[\S8.3.1]{Springer09}.
The claim follows by \eqref{equation: monoid simples} since \(\ell(w)=\#\Phi(w)\).
(c) If \(s\in\tau(w)\) then
\begin{equation*}\label{equation: tau(v*w)}
(v\star w)\star s=v\star(w\star s)=v\star w
\end{equation*}
so \(s\in\tau(v\star w)\).
(d) The relation \(\tau(w_J)=\Phi(w_J)\cap S\) holds by (b).
Observe that \(\Phi(w_J)\) is the set of roots corresponding to the Levi of \(P_J\) by \cite[\S8.4]{Springer09}.
Hence \(\Phi(w_J)\cap S=J\).
\end{proof}
\end{facts}

\begin{proposition}\label{proposition: monoid reduced expressions}
Let \(v=s_1\cdots s_k\) and \(w=t_1\cdots t_\ell\) be reduced expressions.
Then the following hold.
\begin{enumerate}[(a)]
\item \(v=s_1\star\cdots\star s_k\).
\item For some \(1\le i_1<\cdots <i_j\le \ell\),
\begin{equation}
v\star w=s_1\cdots s_k\, t_{i_1}\cdots t_{i_j}
\end{equation}
is a reduced expression.
\item \(\ell(v),\ell(w)\le\ell(v\star w)\le\ell(v)+\ell(w)\).
\item If \(\ell(v\star w)=\ell(v)+\ell(w)\), then \(v\star w=vw\).
\end{enumerate}

\begin{proof}
(a) \(G_v\subseteq\overline{B\dot s_1B\cdots B\dot s_kB}\), so \(v\le s_1\star\cdots\star s_k\).
But \(k=\ell(v)\geq\ell(s_1\star\cdots\star s_k)\) by \eqref{equation: monoid simples}, so \(v=s_1\star\cdots\star s_k\).
(b) This follows from (a) and \eqref{equation: monoid simples}.
(c) This is immediate from (b).
(d) \(vw\leq v\star w\) by Facts~\ref{facts: monoid} part (c).
By (b), the only way \(k+j=k+\ell\) is for every \(1\leq h\leq \ell\), we have \(i_h=h\).
\end{proof}
\end{proposition}

\begin{remark}
If \(\leq_R\) denotes the \emph{right weak order} on \(W\) (as, e.g., in \cite[\S3.1]{Springer09}), then Proposition~\ref{proposition: monoid reduced expressions} tells us that \(v\leq_Rv\star w\), for any \(v,w\in W\).
\end{remark}


\begin{definition}\label{definition: support}
Define the function \(\sigma\) from \(W\) to subsets of \(S\) by
\begin{equation*}\label{equation: support}
\sigma(w)=\set{s\in S\mid s\leq w}
\end{equation*}
called the \emph{support} of \(w\).
\end{definition}

Thus \(s\in\sigma(w)\) if and only if \(G_s\subseteq G_w\).
The `subword property' of Bruhat order \cite[Theorem 2.2.2]{BjornerBrenti05} easily implies that if \(w\in W\), then all reduced expressions for \(w\) contain the same simple reflections.
The support of \(w\) is this set of simple reflections.

\begin{facts}\label{facts: sigma}
The following is a list of simple facts.
\begin{enumerate}[(a)]
\item \(\sigma(w^{-1})=\sigma(w)\).
\item \(\sigma(w_J)=J\).
\item If \(v\le w\) then \(\sigma(v)\subseteq\sigma(w)\).
\item \(\tau(w)\subseteq\sigma(w)\).
\item \(\sigma(v\star w)=\sigma(v)\cup\sigma(w)\).
\end{enumerate}

\begin{proof}
(a) If \(G_s\subseteq G_w\) then \(G_s\subseteq G_{w^{-1}}\), by taking the inverse.
(b) If \(s\) is in \(\sigma(w_J)\) then \(G_s\subseteq G_{w_J}=P_J\).
Hence \(s\) is in \(J\).
If \(s\) is in \(J\) then \(s\le w_J\) so \(s\) is in \(\sigma(w_J)\).
(c) Let \(s\) be in \(\sigma(v)\).
Then \(s\le v\).
Hence \(s\le w\) and \(s\) is in \(\sigma(w)\).
(d) Let \(s\) be in \(\tau(w)\).
Then \(G_wP_s=G_w\) so \(P_s\subseteq G_w\) since \(B\subseteq G_w\).
Hence \(s\) is in \(\sigma(w)\).
(e) Let \(s\) be in \(\sigma(v\star w)\).
Then Proposition~\ref{proposition: monoid reduced expressions} part (b) shows \(s\) is in \(\sigma(v)\) or \(\sigma(w)\).
If \(s\) is in \(\sigma(v)\cup\sigma(w)\) then \(G_s\subseteq G_vG_w=G_{v\star w}\).
Hence \(s\) is in \(\sigma(v\star w)\).
\end{proof}
\end{facts}

The following Lemma will be used later on.

\begin{lemma}\label{lemma: sigma is tau}
Let \(w\) be in \(W\) and \(I\subseteq S\).
The following are equivalent.
\begin{enumerate}[(a)]
\item \(X_w=P_I/B\).
\item \(w=w_I\).
\item \(\sigma(w)=\tau(w)=I\).
\end{enumerate}

\begin{proof}
If \(X_w=P_I/B\) then \(G_w=P_I=G_{w_I}\) so \(w=w_I\).
Suppose that \(w=w_I\) so \(\sigma(w)=I\) and \(\tau(w)=I\) by Facts~\ref{facts: sigma} (b) and Facts~\ref{facts: tau} (d).
Suppose \(\sigma(w)=\tau(w)=I\).
Then \(G_w\subseteq P_I\) since \(G_v\subseteq P_{\sigma(v)}\) always holds.
For every \(s\in I=\tau(w)\), we have \(G_wP_s=G_w\).
Since \(P_I\) is generated by \(\set{P_s\mid s\in I}\) we have \(G_wP_I=G_w\), so \(P_I\subseteq G_w\).
\end{proof}
\end{lemma}

\section{Fiber Bundle Decompositions and BP Decompositions}

Richmond-Slofstra define the notion of BP decomposition in any Coxeter group.
They use this to prove that a necessary and sufficient condition for the morphism \(\pi\,\colon X_w^J\to X_w^I\) to be a fiber bundle, where \(J\subseteq I\), is characterized in terms of BP decompositions.
So, in this case, the geometry of \(X_w^J\) is reduced to the geometry of a Schubert variety in a simpler flag variety \(X_w^I\) and a smaller dimensional \(X_u^J\) in the same flag variety as \(X_w^J\).

In this section, we describe the fiber bundle structure of \(X_w^J\) explicitly as a Bott-Samelson type variety.
We choose to work only with \(X_w\) for notational convenience, but the general case follows directly using, for example, Lemma~\ref{lemma: base change}.
Our main result in this section is Proposition~\ref{proposition: BP isom}, which provides three isomorphisms of the form \(\mu\), as in Proposition~\ref{proposition: mu}, onto \(X_w\) whenever \(w\) admits a BP decomposition.
The various isomorphisms will allow us to: (i) describe \(X_w\) in terms of smaller dimensional Schubert varieties in the same flag variety, (ii) describe \(X_w\) in terms of Schubert varieties in smaller flag varieties, and (iii) force \(\mu\) to satisfy maximal equivariance, as in \eqref{equation: equivariance}.

We apply our perspective to some results from \cite{RichmondSlofstra16} that we will need in later sections.
One main result we will use from \cite[Theorem 3.6]{RichmondSlofstra16} provides a fiber bundle structure for any \(\bQ\)-smooth (also known as rationally smooth) Schubert variety, with base a Schubert variety in a maximal parabolic flag variety.

Using our second isomorphism in Proposition~\ref{proposition: BP isom}, we can iterate this procedure directly to describe their sequence of fiber bundles \cite[Corollary 3.7]{RichmondSlofstra16} as a single Bott-Samelson type variety.
This leads to the definition of complete BP decomposition from \cite{RichmondSlofstra17}, which we use repeatedly in the sequel.


We consider the morphism \(\mu\) of Proposition~\ref{proposition: mu} and give some information on the fiber.
This will be applied to determine when \(\mu\) is birational, and also provide information on \(\tau\) and \(\sigma\).

Let \(u=v\star w\), set \(J=\tau(v)\cap\tau(w^{-1})\), and let \(R=P_J\) be the corresponding standard parabolic subgroup.
By Proposition~\ref{proposition: mu}, the map 
\begin{equation*}\label{equation: m is one mu}
\mu\,\colon G_v\times^RX_w\to X_u,\quad
\mu[g,x]=gx
\end{equation*}
is well-defined and proper.
Note that \(G_vR\subseteq G_vP_{\tau(v)}\subseteq G_v\) by definition of \(\tau(v)\); similarly \(RG_w\subseteq G_w\), so \(R\) stabilizes \(G_v\) and \(G_w\).

\begin{proposition}\label{proposition: m is one fiber}
If \(y\le u\), then 
\begin{equation*}\label{equation: m is one fiber}
\mu^{-1}(\dot yB/B)\isom X_v^J\cap\dot yX_{w^{-1}}^J.
\end{equation*}
\begin{proof}
Let \(\varphi\,\colon G_v\times^RX_w\to X_v^J\times X_u\) by
\begin{equation*}\label{equation: m is one phi}
\varphi[g,x]=(gR/R,gx),
\end{equation*}
as in \eqref{equation: phi}.
The diagram
\[\begin{tikzcd}
G_v\times^RX_w\arrow{rr}{\varphi}\arrow[swap]{dr}{\mu}&&\image(\varphi)\arrow{dl}{\pr_2}\\
&X_u&
\end{tikzcd}\]
commutes.
Therefore,
\begin{equation}\label{equation: fiber isoms}
\begin{split}
\mu^{-1}(\dot yB/B)&\isom\pr_2^{-1}(\dot yB/B)\\
&=\set{(gR/R,\dot yB/B)\in X_v^J\times X_u\mid g^{-1}\dot y\in G_w}\\
&\isom\set{gR/R\in X_v^J\mid g\in \dot yG_{w^{-1}}}\\
&=X_v^J\cap \dot yX_{w^{-1}}^J
\end{split}
\end{equation}
where the second isomorphism is projecting to the first factor.
\end{proof}
\end{proposition}

\begin{corollary}\label{corollary: fiber bundle decomp}
\(\mu\) is an isomorphism if and only if
\begin{equation*}\label{equation: fiber bundle decomp}
\sigma(v)\cap\sigma(w)\subseteq J=\tau(v)\cap\tau(w^{-1}).
\end{equation*}

\begin{proof}
Suppose \(\mu\) is an isomorphism.
We need to show \(\sigma(v)\cap\sigma(w)\subseteq J\).
By \eqref{equation: fiber isoms}, we have \(\mu^{-1}(B/B)\isom X_v^J\cap X_{w^{-1}}^J\), which is equal to the point \(\set{R/R}\) since \(\mu\) is an isomorphism.
Let \(s\in\sigma(v)\cap\sigma(w)\).
Then \(G_s=P_s\subseteq G_v\) by definition of \(\sigma(v)\).
So \(G_sR\subseteq G_vR\).
Also \(G_sR/R\subseteq X_{w^{-1}}^J\) by definition of \(s\in\sigma(w)=\sigma(w^{-1})\).
So \(G_sR/R\subseteq X_v^J\cap X_{w^{-1}}^J=\set{R/R}\), that is, \(G_sR=R\).
So \(s\in \tau(w_J^{-1})=\tau(w_J)=J\).

Conversely, by upper semi-continuity of proper morphisms, it suffices to show that the fiber \(\mu^{-1}(B/B)\) is a point.
Observe \(X_v^J\cap X_{w^{-1}}^J\) is closed and \(B\)-stable, so it is a union of Schubert varieties.
Let \(X_y^J\) be an irreducible component.
But \(\sigma(y)\subseteq\sigma(v)\cap\sigma(w^{-1})\subseteq J\) so \(X_y^J=G_yR/R=R/R\).
Now \(\mu\) is a bijective morphism onto a normal variety, so is an isomorphism by Zariski's Main Theorem.
\end{proof}
\end{corollary}

\begin{remark}\label{remark: fiber bundle decomp}
Recall that \(J=\tau(v)\cap\tau(w^{-1})\subseteq\sigma(v)\cap\sigma(w)\) always holds by Facts~\ref{facts: sigma} (d).
So \(\mu\) is an isomorphism if and only if \(\tau(v)\cap\tau(w^{-1})=\sigma(v)\cap\sigma(w)\).
\end{remark}

\begin{corollary}\label{corollary: push-pull closed}
Let \(I_0,I_1\subseteq S\) and \(w=w_{I_0}\star w_{I_1}\).
For \(R_1=P_{I_0}\cap P_{I_1}\) the standard parabolic subgroup corresponding to \(I_0\cap I_1\), 
\begin{equation*}\label{equation: push-pull closed}
\mu\,\colon P_{I_0}\times^{R_1}P_{I_1}/B\to X_w\
\end{equation*}
is always an isomorphism.
Furthermore, 
\begin{equation}\label{equation: parabolic inflation}
\mu'\,\colon P_{\tau(w^{-1})}\times^RP_{\tau(w)}/B\to X_w,
\end{equation}
where \(R=P_{\tau(w^{-1})}\cap P_{\tau(w)}\), is an isomorphism.
\begin{proof}
The relation \eqref{equation: fiber bundle decomp} holds since \(\sigma(w_I)=\tau(w_I)=\sigma(w_I^{-1})\), so \(\mu\) is an isomorphism.
For the second statement, note that \(w=w_{I_0}\star w_{I_1}\) implies that \(I_0\subseteq\tau(w^{-1})\) and \(I_1\subseteq\tau(w)\), so
\begin{equation*}\label{equation: push-pull closed equivariance}
\begin{split}
w_{\tau(w^{-1})}\star w_{\tau(w)}&=(w_{\tau(w^{-1})}\star w_{I_0})\star(w_{I_1}\star w_{\tau(w)})\\
&=w_{\tau(w^{-1})}\star w\star w_{\tau(w)}\\
&=w.
\end{split}
\end{equation*}
Now the first statement applies to conclude \(\mu'\) is an isomorphism.
\end{proof}
\end{corollary}


In this section, we provide three isomorphisms of the form \(\mu\) to \(X_w\), whenever \(w\) admits a BP decomposition.

Suppose \(J\subseteq S\) and \(w\in W\).
By \cite[Corollary 2.4.5]{BjornerBrenti05}, there exists a unique minimal (with respect to Bruhat order) element \(u_0\) in the coset \(wW_J\).
We may therefore write \(w=u_0u_1\) for \(u_1\) in \(W_J\).
This expression for \(w\) is called the \emph{parabolic decomposition} of \(w\) with respect to \(J\).

\begin{facts}\label{facts: parabolic decomp}
Let \(w=u_0u_1\) be a parabolic decomposition with respect to \(J\).
\begin{enumerate}[(a)]
\item \(\ell(w)=\ell(u_0)+\ell(u_1)\) and \(w=u_0u_1=u_0\star u_1\).
\item \(w=u_0u_1\) is also a parabolic decomposition with respect to \(\sigma(u_1)\).
\item Suppose \(J\subseteq \tau(w)\).
Then \(w=(ww_J^{-1})w_J\) is parabolic with respect to \(J\).
In particular, if \(J\subseteq I\), then \(w_I=(w_Iw_J^{-1})w_J\) is parabolic with respect to \(J\).
\item If \(w=u_0u_1\) is a parabolic decomposition, then \(\tau(u_0)\cap\tau(u_1^{-1})=\emptyset\).
\end{enumerate}
\begin{proof}
(a) See \cite[Proposition 2.4.4]{BjornerBrenti05} for the first statement.
The second statement follows from Proposition~\ref{proposition: monoid reduced expressions}.
(b) The relation \(\sigma(u_1)\subseteq J\) shows that \(u_0\) is also minimal with respect to \(\sigma(u_1)\).
(c) Let \(w=u_0u_1\) be the parabolic decomposition of \(w\) with respect to \(J\).
Suppose (for a contradiction) that \(u_1\) is not equal to \(w_J\).
Then there exists \(s\) in \(J\) such that \(\ell(u_1s)=\ell(u_1)+1\).
Let \(v_0v_1\) be the parabolic decomposition of \(v=ws\) with respect to \(J\).
Then \(v=u_0u_1s\), \(u_0=\min(wsW_J)=\min(wW_J)\), and \(u_1s\) in \(W_J\) force \(v_0=u_0\) and \(v_1=u_1s\) by the uniqueness of parabolic decomposition.
Hence \(\ell(w)=\ell(v_0)+\ell(v_1)=\ell(u_0)+\ell(u_1s)=\ell(u_0)+\ell(u_1)+1\) gives us the desired contradiction, since \(\ell(w)=\ell(u_0)+\ell(u_1)\).
(d) The minimal element \(u_0\) satisfies \(\tau(u_0)\cap J=\emptyset\) since \(\ell(u_0s)>\ell(u_0)\) for every \(s\) in \(J\).
But \(\tau(u_1^{-1})\subseteq\sigma(u_1^{-1})=\sigma(u_1)\subseteq J\) so the claim follows.
\end{proof}
\end{facts}

\begin{lemma}\label{lemma: tau and reduced decompositions}
Let \(u=vw\) be any expression for which \(\ell(u)=\ell(v)+\ell(w)\).
Then \(\tau(u)\subseteq\tau(v)\cup\sigma(w)\).

\begin{proof}
We proceed by induction on \(\ell(w)\).
If \(\ell(w)=0\) then \(u=v\) and there is nothing to prove.

Assume \(\ell(w)\geq 1\) and consider reduced expressions
\begin{equation*}\label{equation: uv reduced}
v=s_1\cdots s_k,\quad w=t_1\cdots t_\ell
\end{equation*}
where \(\ell\geq 1\).
Then \(u=s_1\cdots s_kt_1\cdots t_\ell\) is a reduced expression.
Let \(s\in\tau(u)\smallsetminus\sigma(w)\) and write \(t=t_\ell\).
Claim: \(s\in\tau(ut)\).

Consider \(I=\set{s,t}\) and let \(u=u_0u_1\) be the parabolic decomposition of \(u\) with respect to \(I\).
We have \(s,t\in\tau(u)\), so \(u_1=w_I\) by Facts~\ref{facts: parabolic decomp} (c).
Any reduced expression of \(w_I\) alternates \(s\) and \(t\), so we can find \(y<w\) such that \(w_I=yst\), where \(\ell(w_I)=\ell(y)+\ell(s)+\ell(t)\).
Then \(ut=u_0ys\) such that \(\ell(ut)=\ell(u_0)+\ell(y)+1\), so \(ut=u_0\star y\star s\) by Proposition~\ref{proposition: monoid reduced expressions}.
It follows that \(ut\star s=ut\), i.e., \(s\in\tau(ut)\) as claimed.

Apply the induction hypothesis to \(ut=v(wt)\), so we have \(s\in\tau(v)\) since \(s\notin\sigma(wt)\).
It follows that \(\tau(u)\subseteq\tau(v)\cup\sigma(w)\).
\end{proof}
\end{lemma}

\begin{definition}\label{definition: BP decomp}
A \emph{Billey-Postnikov decomposition} (or \emph{BP decomposition}) of \(w\) with respect to \(I\) is a parabolic decomposition \(w=u_0u_1\) with respect to \(I\) such that
\begin{equation}\label{equation: BP decomp}
\sigma(u_0)\cap I\subseteq\tau(u_1^{-1}).
\end{equation}
Note that a BP decomposition with respect to \(I\) is also a BP decomposition with respect to \(\sigma(u_1)\) since
\begin{equation*}\label{equation: BP decomp'}
\sigma(u_0)\cap\sigma(u_1)\subseteq\tau(u_1^{-1})
\end{equation*}
holds true.
\end{definition}

\begin{proposition}\label{proposition: BP isom}
Suppose \(w=u_0u_1\) is a BP decomposition with respect to some \(I\).
Let \(J=\sigma(u_0)\cap\sigma(u_1)\), \(J'=\tau(u_1^{-1})\), and \(J''=\tau(w_0)\), where \(w_0=u_0\star w_{J'}\).
Let \(v_0=u_0\star w_J\), \(w_1=w_{J'}\star u_1\), and let \(R,R',R''\) be the standard parabolic subgroups corresponding to \(J,J',J''\) respectively.
Then the following hold.
\begin{enumerate}[(i)]
\item
\(w=v_0\star u_1\), \(J=\tau(v_0)\cap\tau(u_1^{-1})\), and
\(\mu\,\colon G_{v_0}\times^RX_{u_1}\to X_w\)
is an isomorphism.
\item
\label{left equivariance}
\(w=w_0\star u_1\), \(J'=\tau(w_0)\cap\tau(u_1^{-1})\), and
\(\mu'\,\colon G_{w_0}\times^{R'}X_{u_1}\to X_w\)
is an isomorphism such that \(\tau(w^{-1})=\tau(w_0^{-1})\).
\item
\label{left then right equivariance}
\(w=w_0\star w_1\), \(J''=\tau(w_0)\cap\tau(w_1^{-1})\), and
\(\mu''\,\colon G_{w_0}\times^{R''}G_{w_1}\to G_w\)
is an isomorphism such that \(\tau(w^{-1})=\tau(w_0^{-1})\) and \(\tau(w)=\tau(w_1)\).
\end{enumerate}

\begin{proof}
(i) It is enough to show that
\begin{equation*}
\sigma(v_0)\cap\sigma(u_1)\subseteq J
\end{equation*}
by \eqref{equation: fiber bundle decomp}.
But \(\sigma(v_0)=\sigma(u_0\star w_J)=\sigma(u_0)\cup J\) by Facts~\ref{facts: sigma}.

(ii) The relation
\begin{equation*}
\sigma(w_0)\cap\sigma(u_1)\subseteq J'
\end{equation*}
holds since \(\sigma(w_0)\cap\sigma(u_1)=\sigma(u_0\star w_{J'})\cap\sigma(u_1)=(\sigma(u_0)\cup J')\cap\sigma(u_1)\) and \(J\subseteq J'\).

By Facts~\ref{facts: tau}, we have \(\tau(w_0^{-1})\subseteq\tau(w^{-1})\) since \(w=w_0\star u_1\).
The other inclusion takes more work.
We prove
\(\tilde\mu'\,\colon G_{w_0'}\times^{R'}X_{u_1}\to X_w\)
is an isomorphism, where \(w_0'=w_{\tau(w^{-1})}\star u_0\star w_{J'}\).
Once this is done, we will have \(w_0=u_0\star w_{J'}\le w_0'\).
But the dimensions of the fiber bundle gives \(\ell(w)=\dim(X_w)=\ell(w_0')-\ell(w_{J'})+\ell(u_1)\) and \(\ell(w)=\dim(X_w)=\ell(w_0)+\ell(u_1)-\ell(w_{J'})\), from \(\mu'\).
Therefore, \(\ell(w_0)=\ell(w_0')\), so \(w_0=w_0'\).
That is, \(w_0=w_{\tau(w^{-1})}\star u_0\star w_{J'}\).
So \(\tau(w^{-1})=\tau(w_{\tau(w^{-1})}^{-1})\subseteq\tau(w_0^{-1})\), and we will have the final statement of (ii).

We now prove that \(\tilde\mu'\) is an isomorphism.
First \(w_0'\star u_1=(w_{\tau(w^{-1})}\star u_0\star w_{J'})\star u_1=w_{\tau(w^{-1})}\star u_0\star (w_{J'}\star u_1)=w_{\tau(w^{-1})}\star(u_0\star u_1)=w_{\tau(w^{-1})}\star w=w\).
Next, \(\tau(w_0')\cap\tau(u_1^{-1})=\tau(u_1^{-1})=J'\), since \(\tau(u_1^{-1})=J'\subseteq\tau(w_0')\) by Facts~\ref{facts: tau}.
Then the condition for isomorphism of \eqref{equation: fiber bundle decomp} 
\begin{equation*}
\begin{split}
\sigma(w_0')\cap\sigma(u_1)&=(\tau(w^{-1})\cup\sigma(u_0)\cup J')\cap\sigma(u_1)\\
&=(\tau(w^{-1})\cap\sigma(u_1))\cup(\sigma(u_0)\cap\sigma(u_1))\cup(J'\cap\sigma(u_1))\\
&\subseteq(\tau(w^{-1})\cap\sigma(u_1))\cup J'\quad\text{ by }\eqref{equation: BP decomp}\\
&\subseteq((\tau(u_1^{-1})\cup\sigma(u_0^{-1}))\cap\sigma(u_1))\cup J'\quad\text{ by }\mathrm{Lemma}~\ref{lemma: tau and reduced decompositions}\\
&\subseteq(J'\cap\sigma(u_1))\cup(\sigma(u_0)\cap\sigma(u_1))\cup J'\\
&\subseteq(J'\cap\sigma(u_1))\cup J'\quad\text{ by }\eqref{equation: BP decomp}\\
&\subseteq J'
\end{split}
\end{equation*}
is satisfied.

(iii) The relation
\begin{equation*}
\sigma(w_0)\cap\sigma(w_1)\subseteq J''
\end{equation*}
holds since \(\sigma(w_0)\cap\sigma(w_1)=(\sigma(u_0)\cup J')\cap(J'\cup \sigma(u_1))\) and \(J\subseteq J'\subseteq J''\).

By Facts~\ref{facts: tau}, we have \(\tau(w_1)\subseteq\tau(w)\) since \(w=w_0\star w_1\).
We prove
\(\tilde\mu''\,\colon G_{w_0}\times^{R''}G_{w_1'}\to G_w\)
is an isomorphism, where \(w_1'=w_1\star w_{\tau(w)}\).
Similar to the proof of (ii), this will give the final statement of (iii).

The subset of simple reflections
\begin{align}
\begin{split}
\sigma(w_0)\cap\sigma(w_1')&=(\sigma(u_0)\cup J')\cap(J'\cup \sigma(u_1)\cup\tau(w))\\
&\subseteq J''\cup ((\sigma(u_0)\cup J')\cap\tau(w))\\
&=J''\cup(\sigma(u_0)\cap\tau(w))\cup(J'\cap\tau(w))
\end{split}
\end{align}
is contained in \(J''\) if and only if \(\sigma(u_0)\cap\tau(w)\subseteq J''\).
Let \(s\) be in \(\sigma(u_0)\cap\tau(w)\).
If \(s\) is in \(\sigma(u_1)\) then \(s\) is in \(\sigma(u_0)\cap\sigma(u_1)=J\subseteq J'\subseteq J''\).
It remains to show that if \(s\) is in \(\sigma(u_0)\cap\tau(w)\) but not \(\sigma(u_1)\) then \(s\) is in \(J''\).

Let \(w_0=u_0w_{J'}\) and \(u_1^{-1}=v_1w_{J'}\) be parabolic decompositions with respect to \(J'\).
Then \(w=u_0u_1=(u_0w_{J'})v_1^{-1}\) satisfies \(\ell(w)=\ell(u_0w_{J'})+\ell(v_1^{-1})\) since \(\ell(u_0)+\ell(u_1)=\ell(u_0)+\ell(w_{J'})+\ell(v_1^{-1})=\ell(u_0w_{J'})+\ell(v_1^{-1})\).
Applying Lemma~\ref{lemma: tau and reduced decompositions} to the above relation shows that \(s\) is in \(\tau(u_0w_{J'})=\tau(u_0\star w_{J'})=\tau(w_0)\) since \(\sigma(v_1^{-1})\subseteq\sigma(u_1)\) by Facts~\ref{facts: sigma}.
But \(J''=\tau(w_0)\) so the claim follows.
\end{proof}
\end{proposition}

\begin{remark}
The first isomorphism to \(X_w\) in Proposition~\ref{proposition: BP isom} (i) is best for providing small dimensional \(X_{v_0}\) and \(X_{u_1}\) in the same flag variety as \(X_w\) (i.e., this is the best chance of giving \(\ell(v_0)<\ell(w)\)).
In this paper, we will most often use the second isomorphism in Proposition~\ref{proposition: BP isom} (ii) because it is best suited for describing \(X_w\) in terms of a Schubert variety in a smaller flag variety \(X_{w_0}^{J'}\) and a smaller dimensional Schubert variety \(X_{u_1}\) in the same flag variety as \(X_w\).
We will use an isomorphism similar to Proposition~\ref{proposition: BP isom} (iii) in the sequel when we need to satisfy \eqref{equation: equivariance}.
However, ensuring that \(\ell(w_1)<\ell(w)\) can require additional care.
\end{remark}


In this section, we recall grassmannian BP decompositions from \cite{RichmondSlofstra16}.
We use Proposition~\ref{proposition: BP isom} to describe all \(\bQ\)-smooth Schubert varieties as Bott-Samelson type varieties.
This leads naturally to the notion of a complete BP decomposition, which is an iterated version of grassmannian BP decomposition.

We recall terminology from \cite{RichmondSlofstra16}.
A \emph{generalized grassmannian} is a flag variety \(X^I\) such that \(\#(S\smallsetminus I)=1\).
A \emph{grassmannian Schubert variety} is a Schubert variety \(X_w^I\) in a generalized grassmannian.
A \emph{grassmannian parabolic decomposition} is a parabolic decomposition of \(w\) with respect to \(I\) such that \(\#(\sigma(w)\cap I)=\#\sigma(w)-1\).
A \emph{grassmannian BP decomposition} is a BP decomposition that is a grassmannian parabolic decomposition.

\begin{facts}\label{facts: grassmannian}
We list some facts describing the terminology, along with some easy facts we will use.
\begin{enumerate}[(a)]
\item
If \(w=u_0u_1\) is a grassmannian parabolic decomposition with respect to \(I\), then \(X_w^I\) is isomorphic to a grassmannian Schubert variety (possibly for a smaller group).
\item
If \(\#(\sigma(w)\smallsetminus\tau(w))\le1\), then \(X_w^{\tau(w)}\) is isomorphic to a grassmannian Schubert variety such that \(X_w\to X_w^{\tau(w)}\) is the base change of the fiber bundle \(X\to X^{\tau(w)}\) with respect to inclusion.
\item
A grassmannian BP decomposition \(w=u_0u_1\) with respect to \(I\) is also a grassmannian BP decomposition of \(w\) with respect to \(\sigma(u_1)\).
\item
If \(w=u_0u_1\) is a grassmannian BP decomposition with respect to \(I\) then \(\#\tau(u_0)=1\).
\end{enumerate}

\begin{proof}
(a) Let \(L_{\sigma(w)}\) be the (connected reductive) Levi subgroup of \(P_{\sigma(w)}\).
Then \(X_w^I\subseteq G/P_I\) is isomorphic to \(X_w^{\sigma(w)\cap I}\subseteq L_{\sigma(w)}/(L_{\sigma(w)}\cap P_{\sigma(w)\cap I})\) since the inclusion of flag varieties gives a closed embedding of Schubert varieties of the same dimension.
Note that \(L_{\sigma(w)}\cap P_{\sigma(w)\cap I}\) is a parabolic subgroup of \(L_{\sigma(w)}\) containing the Borel subgroup \(L_{\sigma(w)}\cap B\), and it corresponds to the simple reflections \(\sigma(w)\cap I\) (e.g., by \cite[\S8.4]{Springer09}).
The flag variety corresponding to the Levi is a generalized grassmannian since \(\#(\sigma(w)\smallsetminus(\sigma(w)\cap I))=1\) by definition of grassmannian parabolic decomposition.

(b) Suppose that \(\#(\sigma(w)\smallsetminus\tau(w))\le1\).
If \(\sigma(w)=\tau(w)\) then \(X_w=P_{\tau(w)}/B\) and \(P_{\tau(w)}/P_{\tau(w)}\) is isomorphic to the minimum Schubert variety in any (grassmannian) flag variety.
Note that \(X_w\to X_w^{\tau(w)}\) is always the base change of the fiber bundle \(X\to X^{\tau(w)}\) with respect to inclusion.

If \(\#(\sigma(w)\smallsetminus\tau(w))=1\) then \(X_w^{\tau(w)}\) is isomorphic to a grassmannian Schubert variety by (a), since the parabolic decomposition of \(w\) with respect to \(\tau(w)\) is grassmannian.

(c) Suppose that \(w=u_0u_1\) is a grassmannian BP decomposition with respect to \(I\).
Then it is a BP decomposition with respect to \(\sigma(u_1)\) by \eqref{equation: BP decomp'}.
So it is enough to show that \(\#(\sigma(w)\smallsetminus\sigma(u_1))=1\) to give the grassmannian condition for parabolic decompositions.
Let \(s\) be the unique simple reflection in \(\sigma(w)\) not in \(I\), by definition of grassmannian parabolic decomposition of \(w\) with respect to \(I\).
Then \(s\) is not in \(\sigma(u_1)\subseteq I\), so \(s\) is in \(\sigma(u_0)\) since \(\sigma(w)=\sigma(u_0)\cup\sigma(u_1)\).
Let \(t\) be any element of \(\sigma(u_0)\) such that \(s\neq t\).
Then \(t\) is in \(I\) by the uniqueness of \(s\).
By definition of BP decomposition, \(\sigma(u_0)\cap I\subseteq\tau(u_1^{-1})\) so \(t\) is in \(\sigma(u_1)\) and the claim follows.

(d) If \(w=u_0u_1\) is a grassmannian BP decomposition with respect to \(I\) then it is with respect to \(\sigma(u_1)\).
Then \(\#(\sigma(w)\smallsetminus\sigma(u_1))=1\) by definition of grassmannian parabolic decomposition.
But \(\tau(u_0)\cap\sigma(u_1)=\emptyset\) by definition of parabolic decomposition.
The claim follows.
\end{proof}

\end{facts}

Richmond-Slofstra \cite{RichmondSlofstra16} show that any \(\bQ\)-smooth Schubert variety \(X_w\) yields a grassmannian BP decomposition \(w=u_0u_1\).
Thus any \(\bQ\)-smooth Schubert variety \(X_w\) is a fiber bundle with base a grassmannian Schubert variety \(X_{u_0}^I\) and fiber \(X_{u_1}\) a smaller \(\bQ\)-smooth Schubert variety.
It follows that the procedure can be applied recursively to reduce the geometry of every \(\bQ\)-smooth Schubert variety to grassmannian Schubert varieties.
We use Proposition~\ref{proposition: BP isom} to describe a resulting Bott-Samelson type structure on every \(\bQ\)-smooth Schubert variety.
Here it is essential that we use Proposition~\ref{proposition: BP isom} (ii) to give us \eqref{equation: a}, and enable a recursive procedure.

\begin{theorem}\label{theorem: Q-smooth isom}
Let \(X_w\) be a \(\bQ\)-smooth Schubert variety.
Then there exists an isomorphism \(\mu\,\colon G_{w_0}\times^{R_1}\cdots\times^{R_m}G_{w_m}/B\to X_w\), such that for every \(0\leq i\leq m\), \(\#\tau(w_i)\geq \#\sigma(w_i)-1\).
We also have \(\tau(w^{-1})=\tau(w_0^{-1})\).

\begin{proof}
The following proof leads to the definition of complete BP decomposition, but could be simplified slightly without this goal in mind.
Suppose \(X_w\) is \(\bQ\)-smooth.
If \(\#\sigma(w)\leq 1\) then \(G_w=P_{\sigma(w)}\) and the theorem is trivial by letting \(\mu\) be the identity map, so assume that \(\#\sigma(w)\geq 2\).

Let \(w=u_0u_1\) be a grassmannian BP decomposition with respect to \(I\), such that \(\#(\sigma(w)\cap I)=\#\sigma(w)-1\), as we can do by \cite[Theorem 3.6]{RichmondSlofstra16}.
Then \(\#\sigma(u_1)\geq 1\) by Facts~\ref{facts: grassmannian} (c) along with the definition of grassmannian BP decomposition.
More preciesly, we have \(\sigma(w)=\set{s}\cup\sigma(u_1)\), where \(\set{s}=\tau(u_0)\) by Facts~\ref{facts: grassmannian} (d).
We also have \(\ell(u_1)<\ell(w)\) since \(\ell(w)=\ell(u_0)+\ell(u_1)\) and \(u_0\neq e\).

Let \(J=\tau(u_1^{-1})\) and \(w_0=u_0\star w_J\), so Proposition~\ref{proposition: BP isom} (ii) gives the isomorphism \(\mu_1\,\colon G_{w_0}\times^{R_1}X_{u_1}\to X_w\) such that \(\tau(w^{-1})=\tau(w_0^{-1})\). 
Since \(X_w\) is \(\bQ\)-smooth and \(\mu_1\) is an isomorphism, \(X_{u_1}\) is \(\bQ\)-smooth.

If \(\#\sigma(u_1)=1\) then \(\#\sigma(u_1)=\#\tau(u_1)\) is a simple reflection and we are done, so assume that \(\#\sigma(u_1)\geq 2\).
Let \(\mu_2\,\colon G_{w_1}\times^{R_2}X_{u_2}\to X_{u_1}\) be an isomorphism such that \(\tau(u_1^{-1})=\tau(w_1^{-1})\) by the above discussion.
Then
\begin{equation*}
\tau(w_0)\cap\tau(u_1^{-1})=\tau(w_0)\cap\tau(w_1^{-1})
\end{equation*}
shows that \(\mu'\,\colon G_{w_0}\times^{R_1}G_{w_1}\times^{R_2}X_{u_2}\to X_w\) is a well-defined morphism.
Hence we argue recursively to get the desired isomorphism \(\mu\).
\end{proof}
\end{theorem}

\begin{definition}\label{definition: complete BP decomp}
A \emph{complete BP decomposition} of \(w\) is a factorization in the Weyl group \(w=u_0\cdots u_m\), where for every \(0\leq i\leq m\), the product \(u_i(u_{i+1}\cdots u_{m+1})\) (with \(u_{m+1}=e\)) is a BP decomposition with respect to \(\sigma(u_{i+1}\cdots u_{m+1})\) such that \(\#\sigma(u_i\cdots u_m)=m+1-i\).
Our definition is equivalent to the original definition in \cite{RichmondSlofstra17} and the definition provided by \cite{AllandRichmond18}.
\end{definition}

In this case, for every \(0\leq i\leq m\), let
\begin{equation*}
\sigma(u_i\cdots u_m)=\set{s_i,\ldots,s_m},
\end{equation*}
where \(s_i\) is the unique simple reflection in \(\sigma(u_i\cdots u_m)\smallsetminus\sigma(u_{i+1}\cdots u_m)\).
For every \(1\leq i\leq m\), let \(J_i=\tau((u_i\cdots u_m)^{-1})\) and set
\begin{equation}\label{equation: w_i}
\begin{split}
w_0&=u_0\star w_{J_1}\\
w_1&=u_1\star w_{J_2}\\
&\,\ \vdots\\
w_{m-1}&=u_{m-1}\star w_{J_m}\\
w_m&=u_m.
\end{split}
\end{equation}

\begin{facts}\label{facts: complete BP decomp}
Let \(\tilde w=(u_0,\ldots,u_m)\) be a complete BP decomposition of \(w\).
\begin{enumerate}[(a)]
\item \(w=w_0\star\cdots\star w_m\).
\item For every \(0\leq i< m\), \(J_{i+1}\subseteq\tau(w_i)\).
\item For every \(1\leq i\leq m\), \(J_i=\tau(w_i^{-1})\subseteq\tau(w_{i-1})\).
We also have \(\tau(w^{-1})=\tau(w_0^{-1})\).
\item For every \(0\leq i\leq m\), \(\tau(w_i)=\sigma(w_i)\) or \(\tau(w_i)=\sigma(w_i)\smallsetminus\set{s_i}\).
\end{enumerate}
\begin{proof}
We prove (d), since the remaining statements follow from above proofs.
First note that for \(0\leq i< m\), \(\sigma(w_i)=\sigma(u_i)\cup\sigma(w_{J_{i+1}})=\sigma(u_i)\cup J_{i+1}=\set{s_i}\cup J_{i+1}\) since \(\sigma(u_i)\cap\sigma(u_{i+1}\cdots u_m)\subseteq\tau((u_{i+1}\cdots u_m)^{-1})=J_{i+1}\) by definition of BP decomposition.
But \(J_{i+1}\subseteq\tau(w_i)\) by (b), so \(\sigma(w_i)=\set{s_i}\cup J_{i+1}\subseteq\set{s_i}\cup\tau(w_i)\) gives the desired statement.
\end{proof}
\end{facts}

\begin{corollary}\label{corollary: complete BP isom}
The map \(\mu\,\colon G_{w_0}\times^{R_1}\cdots\times^{R_m}G_{w_m}/B\to X_w\) is an isomorphism.
\end{corollary}

\begin{remark}
In particular, if a Schubert variety \(X_w\) is \(\bQ\)-smooth, then there exists a complete BP decomposition \(\tilde w=(u_0,\ldots,u_m)\) such that the isomorphism in Theorem~\ref{theorem: Q-smooth isom} is given by Corollary~\ref{corollary: complete BP isom}.
\end{remark}

\section{Small Resolutions}\label{section: small resolutions}

We recall the definition of small resolution from \cite{GoreskyMacPherson83} and we recall a result from \cite{SankaranVanchinathan94} which allows us to change base of a small resolution.
Then we show how a small resolution of the form \(\mu\) for \(X_w\) provides a small resolution for \(X_{w^{-1}}\).
We conclude this section by showing how to glue together small resolutions of the form \(\mu\) to construct new small resolutions.


\begin{definition}
Let \(\widetilde Y\) and \(Y\) be irreducible complex algebraic varieties.
A \emph{resolution of singularities} of \(Y\) is an algebraic morphism \(\xi\,\colon\widetilde Y\to Y\) such that properties (1)-(3) hold true:
(1) \(\xi\) is proper,
(2) \(\xi\) is birational, and
(3) \(\widetilde Y\) is smooth.
A resolution is often required to satisfy: 
(4) \(\xi\) is an isomorphism over the smooth locus of \(Y\),
in which we call it a \emph{strict} resolution of singularities.
\end{definition}

\begin{definition}
A resolution of singularities \(\xi\,\colon\widetilde Y\to Y\) is \emph{small} means for every \(r>0\),
\begin{equation}\label{equation: small}
\codim_Y\set{y\in Y\mid\dim(\xi^{-1}(y))\geq r}>2r,
\end{equation}
where \(\codim_Y(\emptyset)=\infty\).
\end{definition}

A small resolution of a Schubert variety is strict, and can sometimes be used to compute the singular locus (as in \cite{SankaranVanchinathan94}).

It is often easier to describe resolutions in \(G/P\) for \(P\) a maximal parabolic subgroup than to work directly with \(G/B\).
It is then possible to describe explicitly a resolution \(G/B\).
The following appears in a similar form in Sankaran-Vanchinathan \cite[Theorem 2.4]{SankaranVanchinathan94}.

\begin{lemma}{\cite[Theorem 2.4]{SankaranVanchinathan94}}\label{lemma: base change}
Let \(\xi\,\colon\widetilde Y\to Y\) be an algebraic morphism between irreducible varieties and let \(\zeta\,\colon Z\to Y\) be a Zariski locally trivial fiber bundle with irreducible fiber \(F\).
Then base change 
\[
\begin{tikzcd}
\widetilde Y\underset{Y}{\times}Z\arrow{r}{\xi'}\arrow[swap]{d}{\zeta'}&Z\arrow{d}{\zeta}\\
\widetilde Y\arrow{r}{\xi}&Y
\end{tikzcd}
\]
satisfies the following properties.
\begin{enumerate}[(i)]
\item The morphism \(\zeta'\) is a Zariski locally trivial fiber bundle with fiber \(F\).
\item If \(\xi\) is a proper birational algebraic morphism then \(\xi'\) is a proper birational algebraic morphism.
\item Suppose that \(\xi\) and \(\xi'\) are resolutions.
Then \(\xi\) is a small resolution if and only if \(\xi'\) is a small resolution.
\end{enumerate}
\end{lemma}


The Schubert variety \(X_w\) is smooth (or \(\bQ\)-smooth) if and only if \(X_{w^{-1}}\) is smooth (respectively, \(\bQ\)-smooth), but \(X_w\) is not necessarily isomorphic to \(X_{w^{-1}}\), as shown in \cite{RichmondSlofstra16}.
We show that \(X_w\) has a small resolution of the form \(\mu\) if and only if \(X_{w^{-1}}\) has a small resolution of the form \(\mu\).

\begin{lemma}\label{lemma: base change to G}
Let \(\mu\,\colon G_{w_0}\times^{R_1}\cdots\times^{R_m}G_{w_m}/B\to X_w\) and \(\mu'\,\colon G_{w_0}\times^{R_1}\cdots\times^{R_m}G_{w_m}\to G_w\) be given by multiplication.
Then the diagram
\begin{equation*}
\begin{tikzcd}
G_{w_0}\times^{R_1}\cdots\times^{R_m}G_{w_m}\arrow{r}{\mu'}\arrow[swap]{d}{\pi'}&G_w\arrow{d}{\pi}\\
G_{w_0}\times^{R_1}\cdots\times^{R_m}G_{w_m}/B\arrow{r}{\mu}&X_w
\end{tikzcd}
\end{equation*}
is a base change.

\begin{proof}
Let \(Z=G_{w_0}\times^{R_1}\cdots\times^{R_m}G_{w_m}/B\underset{X_w}{\times}G_w\) be the fibered product of \(\mu\) and \(\pi\).
Then the universal property of fibered product provides a morphism \(\alpha\,\colon G_{w_0}\times^{R_1}\cdots\times^{R_m}G_{w_m}\to Z\).
Explicitly, we have
\begin{equation*}
\alpha[g_0,\ldots,g_m]=([g_0,\ldots,g_mB/B],g_0\cdots g_m).
\end{equation*}
Define \(\beta\,\colon Z\to G_{w_0}\times^{R_1}\cdots\times^{R_m}G_{w_m}\) by
\begin{equation*}
\beta([g_0,\ldots,g_{m-1},g_mB/B],g)=[g_0,\ldots,g_{m-1},(g_0\cdots g_{m-1})^{-1}g]
\end{equation*}
which is the morphism induced by quotienting \(\pr_1\,\colon G_{w_0}\times^{R_1}\cdots\times^{R_m}G_{w_m}\underset{X_w}{\times}G_w\to G_{w_0}\times^{R_1}\cdots\times^{R_m}G_{w_m}\).
Then \(\alpha\) and \(\beta\) are inverse algebraic morphisms.
\end{proof}
\end{lemma}

\begin{proposition}\label{proposition: inverse trick}
Let \(\mu\,\colon G_{w_0}\times^{R_1}\cdots\times^{R_m}G_{w_m}/B\to X_w\) be a small resolution of \(X_w\).
Then \(\nu\,\colon G_{w_m^{-1}}\times^{R_m}\cdots\times^{R_1}G_{w_0^{-1}}/B\to X_{w^{-1}}\) is a small resolution of \(X_{w^{-1}}\).

\begin{proof}
Consider the base change
\[
\begin{tikzcd}
G_{w_0}\times^{R_1}\cdots\times^{R_m}G_{w_m}\arrow{r}{\mu'}\arrow{d}{}&G_w\arrow{d}{}\\
G_{w_0}\times^{R_1}\cdots\times^{R_m}G_{w_m}/B\arrow{r}{\mu}&X_w
\end{tikzcd}
\]
by Lemma~\ref{lemma: base change to G}, and consider a similar diagram for \(\nu\).
For \(v\in W\), define \(\alpha_v\,\colon G_v\to G_{v^{-1}}\) by
\begin{equation*}\label{equation: inversion}
\alpha_v(g)=g^{-1}
\end{equation*}
so \(\alpha_v\) is an isomorphism.
Let \(\beta\,\colon G_{w_0}\times^{R_1}\cdots\times^{R_m}G_{w_m}\to G_{w_m^{-1}}\times^{R_m}\cdots\times^{R_1}G_{w_0^{-1}}\) be the map on quotients induced by the various \(\alpha_{w_i}\) and reversing coordinates.
We have a commuting diagram
\begin{equation}\label{equation: inverse diagram in G}
\begin{tikzcd}
G_{w_m^{-1}}\times^{R_m}\cdots\times^{R_1}G_{w_0^{-1}}\arrow{r}{\nu'}\arrow{d}{\beta^{-1}}&G_{w^{-1}}\arrow{d}{\alpha_w^{-1}}\\
G_{w_0}\times^{R_1}\cdots\times^{R_m}G_{w_m}\arrow{r}{\mu'}&G_w
\end{tikzcd}
\end{equation}
so \(\mu'\) is a small resolution if and only \(\nu'\) is a small resolution.
The claim follows by Lemma~\ref{lemma: base change}.
\end{proof}
\end{proposition}

\begin{remark}
Proposition~\ref{proposition: inverse trick} shows that \(X_w\) has a fiber bundle decomposition if and only if \(X_{w^{-1}}\) has a fiber bundle decomposition, since \eqref{equation: inverse diagram in G} shows \(\mu\) is an isomorphism if and only if \(\nu\) is an isomorphism (regardless of all \(G_{w_i}\) being smooth).
However, as remarked in \cite{RichmondSlofstra16}, a BP decomposition of \(w\) does not necessarily give a BP decomposition of \(w^{-1}\).
\end{remark}


A fiber bundle decomposition of \(X_w\) allows us to glue small resolutions of the form \(\mu\), if we assume some compatibility with equivariance.

\begin{lemma}\label{lemma: fiber bundle R}
Let \(R\) be any standard parabolic group stabilizing \(G_v\) by right multiplication, and let \(F\) be a left \(R\)-variety.
Then \(\rho\,\colon G_v\times F\to G_v\times^RF\) is a fiber bundle with fiber \(R\).
\end{lemma}

\begin{key lemma}\label{lemma: small glue}
Let \(\nu\,\colon G_{v_0}\times^{R_1'}\cdots\times^{R_n'}G_{v_n}/B\to X_v\) be a small resolution of \(X_v\), and let \(\mu\,\colon G_{w_0}\times^{R_1}\cdots\times^{R_m}G_{w_m}/B\to X_w\) be a small resolution of \(X_w\).
Then \(G_{v_0}\times^{R_1'}\cdots\times^{R_n'}G_{v_n}\times^RG_{w_0}\times^{R_1}\cdots\times^{R_m}G_{w_m}/B\to G_v\times^RX_w\) is a small resolution, where \(R\) is the parabolic subgroup corresponding to \(\tau(v_n)\cap\tau(w_0^{-1})\).

\begin{proof}
Lemmas~\ref{lemma: base change} and \ref{lemma: base change to G} show that \(\nu\) is a small resolution if and only if the base change \(\nu'\) to \(G_v\) is a small resolution.
The diagram
\[
\begin{tikzcd}
G_{v_0}\times^{R_1'}\cdots\times^{R_n'}G_{v_n}\times G_{w_0}\times^{R_1}\cdots\times^{R_m}G_{w_m}/B\arrow{r}{(\nu',\mu)}\arrow{d}{\pi'}&G_v\times X_w\arrow{d}{\pi}\\
G_{v_0}\times^{R_1'}\cdots\times^{R_n'}G_{v_n}\times^R G_{w_0}\times^{R_1}\cdots\times^{R_m}G_{w_m}/B\arrow{r}{[\nu',\mu]}&G_v\times^RX_w
\end{tikzcd}
\]
is a base change of fiber bundles with fiber \(R\).
Since \((\nu',\mu)\) is a small resolution then \([\nu',\mu]\) is a small resolution by Lemma~\ref{lemma: base change}.
\end{proof}
\end{key lemma}

We will apply Lemma~\ref{lemma: small glue} to Schubert varieties of the form \(X_u\isom G_v\times^RX_w\).
Note it is essential to check that small resolutions of \(X_v\) and \(X_w\) give rise to the same \(R\) from Lemma~\ref{lemma: small glue} as in the isomorphism of \(X_u\).
It would be interesting (and useful) to know whether every Schubert variety \(X_w\) admitting a small resolution, admits a (possibly different) small resolution that satisfies \eqref{equation: equivariance}.

\section{Gelfand-MacPherson Resolutions}\label{section: Gelfand-MacPherson}

We recall the construction in \cite[\S2.11]{GelfandMacPherson82} providing a resolution of singularities for any Schubert variety \(X_w^I\). 
The resolution is uniquely determined by subsets of simple reflections \(I_0,\ldots,I_m\) that they call resolution data.
These resolutions are described as iterated base changes of flag varieties, which enables us to compute fibers explicitly.
In particular, \cite{SankaranVanchinathan94} provides a formula for all fiber dimensions.

In this section, we consider Schubert varieties \(X_w^I\) such that
(i) \(w\) is maximal in its \(W_I\)-coset,
(ii) \((W_{\sigma(w)},\sigma(w))\) is a simply laced Coxeter system,
and say that \(X_w^I\) is a \emph{simply laced} Schubert variety.
Note that condition (i) is without loss of generality since \(X_w^I=X_v^I\) whenever \(wW_I=vW_I\).
We show that \emph{any} resolution \(\mu\,\colon G_{w_0}\times^{R_1}\cdots\times^{R_m}X_{w_m}^I\to X_w^I\) is isomorphic to a Gelfand-MacPherson resolution.
In other words, there exists a commuting diagram
\[
\begin{tikzcd}
P_{I_0}\times^{R_1'}\cdots\times^{R_n'}P_{I_n}/P_I\arrow{rr}{\mu'}\arrow[swap]{dr}{\nu}&&G_{w_0}\times^{R_1}\cdots\times^{R_m}X_{w_m}^I\arrow{dl}{\mu}\\
&X_w^I&
\end{tikzcd}
\]
for some \(I_0,\ldots,I_n\) such that \(\mu'\) is an isomorphism.
We further show that if we take \(X_w^I\) smooth and \(\mu\) the identity (resolution) morphism, then \(\mu'\) satisfies \eqref{equation: equivariance}.


\begin{definition}
A sequence \((I_i)\) of sets of simple reflections 
\begin{equation*}\label{equation: resolution data}
\emptyset=I_0,I_1,\ldots,I_m=I
\end{equation*}
is called \emph{resolution data} for the Schubert variety \(X_w^I\) if the iterated base change
\begin{equation*}
Z(I_i):=(P_{I_0}/P_{I_0})\underset{X^{I_0}}{\times}X^{I_0\cap I_1}\underset{X^{I_1}}{\times}\cdots\underset{X^{I_{m-1}}}{\times}X^{I_{m-1}\cap I_m}
\end{equation*}
projects birationally onto \(X_w^I\).
\end{definition}

We remark that the original definition of resolution data is given in terms of the Grothendieck group of a subcategory of a derived category of sheaves on \(X^I\) with the analytic topology, but is equivalent to the definition given here.
When all sets \(I_1,\ldots,I_{m-1}\) have one element, it is the Demazure resolution.

\begin{lemma}\label{lemma: resolution data}
Let \((I_i)\) be resolution data for \(X_w^I\).
Let \(\mu\,\colon P_{I_0}\times^{R_1}P_{I_1}\times^{R_2}\cdots\times^{R_m}P_{I_m}/P_{I_m}\to X^{I_m}\) be given by
\begin{equation*}
\mu[g_0,\ldots,g_mP_{I_m}/P_{I_m}]=g_0\cdots g_mP_{I_m}/P_{I_m}
\end{equation*}
where for every \(1\le i\le m\), \(R_i=P_{I_{i-1}}\cap P_{I_i}\).
Let \(\varphi\,\colon P_{I_0}\times^{R_1}P_{I_1}\times^{R_2}\cdots\times^{R_m}P_{I_m}/P_{I_m}\to Z(I_i)\) be given by
\begin{equation*}
\varphi[g_0,\ldots,g_mP_{I_m}/P_{I_m}]=(P_{I_0}/P_{I_0},g_0R_1/R_1,\ldots,g_0\cdots g_{m-1}R_m/R_m)
\end{equation*}
and \(\nu\,\colon Z(I_i)\to X_w^I\) projection.
Then \(\varphi\) is an isomorphism such that the diagram
\[
\begin{tikzcd}
P_{I_0}\times^{R_1}\cdots\times^{R_m}P_{I_m}/P_{I_m}\arrow[swap]{dr}{\mu}\arrow{rr}{\varphi}&&Z(I_i)\arrow{dl}{\nu}\\
&X^{I_m}&
\end{tikzcd}
\]
commutes.
As a consequence, we have \(w=w_{I_0}\star\cdots\star w_{I_m}\).

\begin{proof}
The proof is the same as that of \eqref{equation: phi}.
\end{proof}
\end{lemma}

\begin{example}\label{example: simple zelevinskii}
We provide a simple example that is one of Zelevinski\u\i's small resolutions.
Let \(W=S_4\) of type \(A_3\).
Consider \(w=\text{( 4 2 3 1 )}\), as in \S\ref{section: A}, so \(\tau(w)=\set{1,3}\).
We have
\begin{equation*}
X_w^{\tau(w)}=\set{E\in\Gr_2(\bC^4)\mid \dim(\bC^2\cap E)\geq 1}.
\end{equation*}
Let \(I_1=\set{1,3}\), \(I_2=\set{2,3}\), \(I_3=\set{1,3}\).
Then \(Z(I_i)\) may be identified with the diagram
\[
\begin{tikzcd}
\bC^4\arrow[dash]{d}{}\\
\bC^3\arrow[dash]{d}{}\\
\bC^2\arrow[dash]{d}{}&F^2\arrow[equal]{l}{}\arrow[dash]{uul}{}&E^2\arrow[dash]{d}{}\arrow[dash]{uull}{}\\
\bC^1\arrow[dash]{d}{}&F^1\arrow[dash]{u}{}\arrow[dash]{dl}{}&E^1\arrow[equal]{l}{}\arrow[dash]{dll}{}\\
0
\end{tikzcd}
\]
where vertical flags are coordinates of \(Z(I_i)\), and horizontal equal signs provide the fibered product relations.
So
\begin{equation*}
Z(I_i)=\set{(F^1,E^2)\in\Gr_1(\bC^2)\times\Gr_2(\bC^4)\mid F^1\subseteq E^2}
\end{equation*}
gives the projection \(\pr_2\,\colon Z(I_i)\to X_w^{\tau(w)}\).

We also have the isomorphism \(\varphi\,\colon P_{I_0}\times^{R_1}P_{I_1}\times^{R_2}P_{I_2}\times^{R_2}P_{I_3}/P_{I_3}\to Z(I_i)\) by (the proof of) Lemma~\ref{lemma: resolution data}.
So \(Z(I_i)\) is a smooth, irreducible, iterated fiber bundle.
Then
\begin{equation*}
\pr_2^{-1}(E)=\begin{cases} \set{\bC^2\cap E},& E\neq\bC^2,\\ \Gr_1(\bC^2),&E=\bC^2,\end{cases}
\end{equation*}
shows \((I_i)\) is resolution data, since the fiber of \(\pr_2\) is generically a point.
Moreover, it shows that \(\pr_2\) is a small resolution.
Note we could (and typically do) use \(\varphi\) to conclude that \((I_i)\) is resolution data by counting dimensions and applying Proposition~\ref{proposition: mu}.
\end{example}


\begin{lemma}\label{lemma: smooth grassmannians}
Let \(X_w\) be a simply laced Schubert variety.
Let \(t\) be a simple reflection in \(\sigma(w)\) and \(I=\sigma(w)-\set{t}\).
Let \(w=u_0u_1\) be the parabolic decomposition of \(w\) with respect to \(I\).
Then \(X_w^I\) is smooth if and only if \(u_0\star w_J=w_{\sigma(u_0)}\), where \(J=I\cap\sigma(u_0)=\sigma(u_0)\smallsetminus\set{t}\).

\begin{proof}
Richmond-Slofstra \cite[Theorem 3.8]{RichmondSlofstra16} show in the simply laced case that \(X_{u_0}^J\) is smooth if and only if \(u_0\) is the maximum element of the minimal length representatives of \(W_{\sigma(u_0)}/W_J\).
By considering maximal length representatives, this is equivalent to \(u_0\star w_J=\max(W_{\sigma(u_0)})=w_{\sigma(u_0)}\) since the function \(w\mapsto u_0\) from \(W\) to minimal length representatives is order preserving by \cite[Proposition 2.5.1]{BjornerBrenti05} (and \(u_0\star w_J\) is always the maximum of \(u_0W_J\)).

If \(X_w^I\) is smooth then \(X_{w\star w_I}\) is smooth since it is the pull-back of \(X_w^I\).
Observe that \(w=u_0\star u_1\) by Facts~\ref{facts: parabolic decomp}, so \(w\star w_I=u_0\star u_1\star w_I=u_0\star w_I=u_0\star w_J\star w_I\) since \(\sigma(u_1)\cup J\subseteq I\).
But \(\mu\,\colon G_{u_0\star w_J}\times^{P_J}P_I/B\to X_{w\star w_I}\) is an isomorphism by \eqref{equation: fiber bundle decomp}, since \(J=\sigma(u_0)\cap I\).
Hence \(X_{u_0}^J\) is smooth, so \(u_0\star w_J=w_{\sigma(u_0)}\) by the previous paragraph.

If \(u_0\star w_J=w_{\sigma(u_0)}\) then the first paragraph shows \(X_{u_0}^J\) is smooth, and the isomorphism \(\mu\) (from the previous paragraph) shows that \(X_{w\star w_I}\) is smooth, and so is \(X_{w\star w_I}^I=X_w^I\).
\end{proof}
\end{lemma}

\begin{theorem}\label{theorem: simply laced smooth isom}
Let \(X_w\) be a simply laced Schubert variety.
\begin{enumerate}[(i)]
\item
\(X_w\) is smooth if and only if there exists an isomorphism \(\mu\,\colon P_{I_0}\times^{R_1}\cdots\times^{R_m}P_{I_m}/B\to X_w\) such that \(\tau(w^{-1})=I_0\).
\item
\(X_w\) is smooth if and only if there exists an isomorphism \(\mu'\,\colon P_{I_0'}\times^{R_1'}\cdots\times^{R_{m'}'}P_{I_{m'}'}\to G_w\) such that \(\tau(w^{-1})=I_0'\) and \(\tau(w)=I_{m'}'\).
\end{enumerate}

\begin{proof}[Proof of (i)]
If there exists such an isomorphism then \(X_w\) is smooth by Proposition~\ref{proposition: mu}.

If \(X_w\) is smooth then \(X_w\) is \(\bQ\)-smooth so we can apply \cite{RichmondSlofstra16} (as in Theorem~\ref{theorem: Q-smooth isom}) to get a complete BP decomposition \(\tilde w=(u_0,\ldots,u_n)\). 
Recall the definition of \(s_0,\ldots,s_n\) and \(w_0,\ldots,w_n\) from Facts~\ref{facts: complete BP decomp}, and for every \(0\le i\le n\), let \(J_i=\tau((u_i\cdots u_n)^{-1})\).
The isomorphism \(\mu\,\colon G_{w_0}\times^{P_{J_1}}\cdots\times^{P_{J_n}}X_{w_n}\to X_w\) from Corollary~\ref{corollary: complete BP isom} is such that for every \(0\le i\le n\), we have \(J_i=\tau(w_i^{-1})\) and \(\#\tau(w_i)\geq \#\sigma(w_i)-1\), by Facts~\ref{facts: complete BP decomp} (d).
Note that \(\tau(w^{-1})=\tau(w_0^{-1})\) by Facts~\ref{facts: complete BP decomp} (c).

The smoothness of \(X_w\) along with fiber bundle structures implies that for every \(0\le i\le n\), \(X_{w_i}\) and \(X_{w_i}^{\tau(w_i)}\) are smooth.
By Lemma~\ref{lemma: smooth grassmannians}, for every \(0\le i< n\), we have \(w_i=w_{\sigma(u_i)}\star w_{J_{i+1}}=w_{\tau(w_i^{-1})}\star w_{\tau(w_i)}\) and \(w_n=s_n=w_{\tau(w_n)}=w_{\tau(w_n^{-1})}\star w_{\tau(w_n)}\).
Hence a repeated application of \eqref{equation: parabolic inflation} shows that
\begin{equation*}
G_{w_0}\times^{P_{J_1}}\cdots\times^{P_{J_n}}X_{w_n}\isom P_{\tau(w_0^{-1})}\times^{R_1}P_{\tau(w_0)}\times^{R_2}\cdots\times^{R_{m-1}}P_{\tau(w_n^{-1})}\times^{R_m}P_{\tau(w_n)}/B
\end{equation*}
where all \(R_i\) are intersections of neighboring parabolic subgroups, \(m=2n+1\), and for every \(1\le i\le n\), \(R_{2i}=P_{J_i}\).
\end{proof}

\begin{proof}[Proof of (ii)]
Consider the set
\(A=\set{0\le i\le n\mid w_{\tau(w^{-1})}\star w_i\neq w_{\tau(w^{-1})}}\)
depending on the complete BP decomposition of \(w\).
If \(A\) is empty then for every \(0\le i\le n\), \(w_{\tau(w^{-1})}\star w_i=w_{\tau(w^{-1})}\) forces \(\sigma(w_i)\subseteq\tau(w^{-1})\).
But \(w=w_0\star\cdots\star w_n\) (by Facts~\ref{facts: complete BP decomp}), so 
\begin{equation*}
\sigma(w)=\bigcup_{i=0}^n\sigma(w_i)\subseteq\tau(w^{-1})
\end{equation*}
by Facts~\ref{facts: sigma}. 
Hence \(w=w_{\sigma(w)}\) by Lemma~\ref{lemma: sigma is tau}.
In this case, \(\mu=\id\) is the desired isomorphism of \(X_w\).

From now on,
\begin{equation}\label{equation: A}
A\neq\emptyset
\end{equation}
is a running assumption.
Define \(k=\min(A)\) so by \eqref{equation: A}, \(0\le k\le n\) is well-defined.

We proceed by induction on \(\ell(w)\).
Let \(\ell=\ell(w)\).
Assume for every \(u\) such that \(X_u\) is smooth and \(\ell(u)<\ell\), then there exists an isomorphism \(\mu_u\,\colon P_{I_0^u}\times^{R_1^u}\cdots\times^{R_{n'}^u}P_{I_{n'}^u}/B\to X_u\) such that \(\tau(u^{-1})=I_0^u\) and \(\tau(u)=I_{n'}^u\).

Observe that
\begin{equation*}
\begin{split}
w&=w_{\tau(w^{-1})}\star w\\
&=w_{\tau(w^{-1})}\star w_0\star\cdots\star w_n\\
&=w_{\tau(w^{-1})}\star w_k\star\cdots\star w_n
\end{split}
\end{equation*}
since for every \(0\le i<k\), we have \(w_{\tau(w^{-1})}\star w_i=w_{\tau(w^{-1})}\). 

Set \(I=\tau(w^{-1})\smallsetminus\set{s_k}\) and let \(u=w_I\star (u_{k+1}\cdots u_n)\).
We claim that \(w_{\tau(w^{-1})}\star u=w\).
The claim is equivalent to showing \(w_{\tau(w^{-1})}\star u_{k+1}\star\cdots\star u_n=w\) since \(I\subseteq\tau(w^{-1})\).
The relation \(\sigma(u_k)\subseteq\tau(w_k^{-1})=\tau((u_k\cdots u_n)^{-1})=J_k\) holds since \(w_k=w_{\sigma(u_k)}\star w_{J_{k+1}}\).
Hence \(\sigma(u_k)\subseteq\tau(w_{k-1})\) since \(w_{k-1}=u_{k-1}\star w_{J_k}\).
The definition of \(k\) forces \(w_{\tau(w^{-1})}\star w_{k-1}=w_{\tau(w^{-1})}\) and so \(\sigma(w_{k-1})\subseteq\tau(w^{-1})\).
Hence \(\sigma(u_k)\subseteq\tau(w_{k-1})\subseteq\sigma(w_{k-1})\subseteq\tau(w^{-1})\).
So
\begin{equation}\label{equation: k+1}
\begin{split}
w&=w_{\tau(w^{-1})}\star w\\
&=w_{\tau(w^{-1})}\star w_0\star\cdots\star w_n\\
&=w_{\tau(w^{-1})}\star w_k\star\cdots w_n\\
&=w_{\tau(w^{-1})}\star w_{\sigma(u_k)}\star w_{J_{k+1}}\star w_{k+1}\star\cdots \star w_n\\
&=w_{\tau(w^{-1})}\star w_{k+1}\star\cdots\star w_n\\
&=w_{\tau(w^{-1})}\star u_{k+1}\star\cdots\star u_n\\
&=w_{\tau(w^{-1})}\star u_{k+1}\cdots u_n
\end{split}
\end{equation}
gives the claim.

By Corollary~\ref{corollary: fiber bundle decomp}, the morphism \(\mu'\,\colon P_{\tau(w^{-1})}\times^{P_I}X_u\to X_w\) is an isomorphism since
\begin{equation*}
\begin{split}
\tau(w^{-1})\cap\sigma(u)&=\tau(w^{-1})\cap(I\cup\sigma(u_{k+1}\cdots u_n))\\
&=(\tau(w^{-1})\cap I)\cup(\tau(w^{-1})\cap\sigma(u_{k+1}\cdots u_n))
\end{split}
\end{equation*}
is contained in \(I\).

We see that \(\tau(u)=\tau(w)\) as follows.
The morphism \(\mu''\,\colon P_{\tau(w^{-1})}\times^{P_I}X_{u\star w_{\tau(w)}}\to X_w\) is onto \(X_w\) since \(w_{\tau(w^{-1})}\star u\star w_{\tau(w)}=w\star w_{\tau(w)}=w\).
By Corollary~\ref{corollary: fiber bundle decomp}, \(\mu''\) is an isomorphism if and only if \(\tau(w^{-1})\cap(\sigma(u)\cup\tau(w))\subseteq I\).
But this holds since \(\tau(w^{-1})\cap\tau(w)\subseteq\tau(w^{-1})\smallsetminus\set{s_k}=I\) by Lemma~\ref{lemma: s_k not in tau(w)} below.
By comparing dimensions with \(\mu'\), the relation \(\tau(w)\subseteq\tau(u)\) must hold.
So \(\tau(u)=\tau(w)\) since \(\tau(u)\subseteq\tau(w)\) by Facts~\ref{facts: tau} (c).

Then \(X_u\) is smooth (by the isomorphism \(\mu'\)) such that \(\ell(u)<\ell(w)=\ell\) (since \(s_k\) is in \(\sigma(w)\smallsetminus\sigma(u)\)) and \(\tau(u)=\tau(w)\).
By the induction hypothesis, there exists an isomorphism \(\mu_u\,\colon P_{I_0^u}\times^{R_1^u}\cdots\times^{R_{n'}^u}P_{I_{n'}^u}/B\to X_u\) such that \(\tau(u^{-1})=I_0^u\) and \(\tau(u)=I_{n'}^u\).
Therefore, \(\mu'\) and \(\mu_u\) give an isomorphism \(P_{\tau(w^{-1})}\times^{P_I}P_{I_0^u}\times^{R_1^u}\cdots\times^{R_{n'}^u}P_{I_{n'}^u}/B\to X_w\) such that \(\tau(w)=\tau(u)=I_{n'}^u\).
\end{proof}
\end{theorem}

To complete the proof of Theorem~\ref{theorem: simply laced smooth isom} we need to prove Lemma~\ref{lemma: s_k not in tau(w)}.
For this we first need a definition.

\begin{definition}
Define a function \(\partial\) from \(W\) to subsets of simple reflections
\begin{equation*}
\partial(w)=\set{s\in S\mid s\notin\sigma(w),\ \exists t\in\sigma(w),\ st\neq ts},
\end{equation*}
called the \emph{boundary} of \(w\).
\end{definition}

Then \(\partial(w)\) is the set of simple reflections which are adjacent to \(\sigma(w)\) in the Coxeter graph of \(W\).
Note that \(\partial(w)=\partial(w^{-1})\), for each \(w\in W\).
In the proof of Theorem~\ref{theorem: simply laced smooth isom}, recall that \(X_w\) is a smooth simply laced Schubert variety, \(\tilde w=(u_0,\ldots,u_n)\) is a complete BP decomposition of \(w\), \(A=\set{0\le i\le n\mid w_{\tau(w^{-1})}\star w_i\neq w_{\tau(w^{-1})}}\), and \(k=\min(A)\) when \(A\) is nonempty.

\begin{lemma}\label{lemma: s_k not in tau(w)} 
\(s_k\notin\tau(w)\).

\begin{proof}
Let \(w^{-1}=v_0v_1\) be the parabolic decomposition of \(w^{-1}\) with respect to \(\tau(w^{-1})\).
Then \(v_1=w_{\tau(w^{-1})}\) by Facts~\ref{facts: parabolic decomp}.
Note \(w=v_1^{-1}v_0^{-1}=w_{\tau(w^{-1})}v_0^{-1}\).
We will show that
\begin{equation*}\label{equation: boundary claim}
s_k\in\partial(v_0).
\end{equation*}
Recall that
\begin{equation}
w=w_{\tau(w^{-1})}\star w_{k+1}\star\cdots\star w_n
\end{equation}
by \eqref{equation: k+1}.
Hence Proposition~\ref{proposition: monoid reduced expressions} gives \(\sigma(v_0)\subseteq \sigma(w_{k+1}\star\cdots\star w_n)\), so \(s_k\) is not in \(\sigma(v_0)\).

We are reduced to showing that there exists \(t\) in \(\sigma(v_0)\) such that \(s_kt\neq ts_k\).
Let \((w_{\tau(w^{-1})}\star w_k)^{-1}=v_0'v_1'\) be the parabolic decomposition with respect to \(\tau(w^{-1})\).
Note \(v_1'=w_{\tau(w^{-1})}\) and \(w_{\tau(w^{-1})}\star w_k=w_{\tau(w^{-1})}(v_0')^{-1}\).
The relation \(v_0'\le v_0\) holds by \cite[Proposition 2.5.1]{BjornerBrenti05} since \((w_{\tau(w^{-1})}\star w_k)^{-1}\le w^{-1}\).
Hence \(\sigma(v_0')\subseteq\sigma(v_0)\) and we show that there exists \(t\) in \(\sigma(v_0')\) such that \(s_kt\neq ts_k\).

Note that the Coxeter graph of \(\sigma(u_k)\) is connected since \(\tau(u_k)=\set{s_k}\) is a single element by definition of grassmannian BP decomposition \(u_k(u_{k+1}\cdots u_n)\).
Let \(K\) be the connected component of the Coxeter graph of \(\sigma(w_k)\) such that \(s_k\) is in \(K\).
In particular, \(\sigma(u_k)\subseteq K\).
Let \(s\) be in a connected component of \(\sigma(w_k)\) other than \(K\).
Then \(\sigma(w_k)=\sigma(u_k\star w_{\tau(w_k)})=\sigma(u_k)\cup\tau(w_k)\) shows \(s\) is in \(\tau(w_k)\).
We also have \(s\star w_k=s\star u_k\star w_{\tau(w_k)}=u_k\star s\star w_{\tau(w_k)}=u_k\star w_{\tau(w_k)}\) since \(s\) is not adjacent to \(\sigma(u_k)\).
Hence \(s\) is in \(\tau(w_k^{-1})\).
So \(K\not\subseteq\sigma(u_k)\) since we would have \(\sigma(w_k)\subseteq\tau(w_k^{-1})\) which contradicts \(w_{\tau(w^{-1})}\star w_k\neq w_{\tau(w^{-1})}\) and \(\tau(w_k^{-1})\subseteq\tau(w^{-1})\) by definition of \(k\).

The previous paragraph shows that we can take a path of minimal length \(t_1,\ldots,t_h\) from \(s_k=t_1\) to \(t_h\in\partial(u_k)\) such that for every \(1\le i\le h\), we have \(t_i\) in \(K\).
Note that \(\set{t_2,\ldots,t_h}\subseteq\tau(w_k)\) since \(\tau(w_k)=\sigma(w_k)\smallsetminus\set{s_k}\).
We see that \(t_h\) is not in \(\tau(w^{-1})\) as follows.
Let \(t_h=s_j\) for some \(k<j\le n\), where \(s_j\) is the unique simple reflection in \(\sigma(u_j)\smallsetminus\sigma(u_{j+1}\cdots u_n)\).
Then \(t_h\) is not in \(\tau(w_k^{-1})\) by \cite[Lemma 6.4]{RichmondSlofstra16} since \(u_k(u_{k+1}\cdots u_n)\) is a parabolic decomposition, \(t_h\) is in \(\partial(u_k)\), and \(\tau(w_k^{-1})=\tau((u_{k+1}\cdots u_n)^{-1})\) by Proposition~\ref{proposition: BP isom}.
For every \(0\le i<k\), \(t_h\) is not in \(\sigma(w_i)\) since \(\sigma(w_i)=\sigma(u_i\star w_{\tau(w_{i+1}^{-1})})=\set{s_i}\cup\tau(w_{i+1}^{-1})\) such that \(s_i\neq s_j\) by definition of complete BP decomposition.
It follows that \(t_h\) is not in \(\tau(w^{-1})\) since \(\tau(w^{-1})=\tau(w_0^{-1})\).

Identifying \(W_{\set{t_2,\ldots,t_h}}\) with \(S_{h}\) of type \(A_{h-1}\), permutations of \(\set{1,\ldots,h}\), (by assuming the path is of minimal length in a simply connected Coxeter graph) gives 
\begin{equation*}\label{equation: A_h}
t_2\cdots t_h=(\ 2\ 3\ \cdots\ h\ 1\ )
\end{equation*}
as in \S\ref{section: A}.
So \(\tau(t_2\cdots t_h)=\set{t_h}\) and hence \(t_2\cdots t_h\) is minimal with respect to \(\tau(w^{-1})\).
Therefore the relation \((w_{\tau(w^{-1})}\star(t_h\cdots t_2))^{-1}\le (w_{\tau(w^{-1})}\star w_{\tau(w_k)})^{-1}=(w_{\tau(w^{-1})}\star w_k)^{-1}\) shows that \(t_2\cdots t_h\le v_0'\) by \cite[Proposition 2.5.1]{BjornerBrenti05}. 

Setting \(t=t_2\) gives our claim as follows.
We have \(\sigma(t_2\cdots t_h)\subseteq\sigma(v_0')\) by the end of the last paragraph, and hence \(t\) is in \(\sigma(v_0')\).
But \(t\) is in \(\partial(s_k)=\partial(t_1)\) since \(K\) is connected.
\end{proof}
\end{lemma}
 
\begin{corollary}\label{corollary: simply laced Gelfand-MacPherson resolution}
Let \(X_w^I\) be a simply laced Schubert variety.
Suppose that \(\mu\,\colon G_{w_0}\times^{R_1}\cdots\times^{R_m}X_{w_m}^I\to X_w^I\) is a resolution of singularities.
Then there exists resolution data \((I_i)\) for \(X_w^I\) and an isomorphism \(\varphi\,\colon Z(I_i)\to G_{w_0}\times^{R_1}\cdots\times^{R_m}X_{w_m}^I\) such that the diagram
\[
\begin{tikzcd}
Z(I_i)\arrow{rr}{\varphi}\arrow[swap]{dr}{\nu}&&G_{w_0}\times^{R_1}\cdots\times^{R_m}X_{w_m}^I\arrow{dl}{\mu}\\
&X_w^I&
\end{tikzcd}
\]
commutes.

\begin{proof}
Without loss of generality, assume that \(I\subseteq\tau(w)\).
For every \(0\le i\le m\), \(G_{w_i}\) is smooth since \(\mu\) is a resolution of singularities.
Then \(X_{w_i}\) is simply laced since \(w_i\le w\) by Facts~\ref{facts: sigma} (c), so the claim follows by Theorem~\ref{theorem: simply laced smooth isom} and Theorem~\ref{lemma: small glue}.
\end{proof}
\end{corollary}

\begin{remark}
By \cite[Remark 3.4]{Perrin07}, Corollary~\ref{corollary: simply laced Gelfand-MacPherson resolution} shows that all resolutions constructed in Perrin \cite{Perrin07} are of the form Gelfand-MacPherson, for some resolution data. 
\end{remark}

\begin{example}
\label{counter example C2}
We provide an example to show that if \(X_w\) is not simply laced, the conclusions of Lemma~\ref{lemma: smooth grassmannians}, Theorem~\ref{theorem: simply laced smooth isom}, and Corollary~\ref{corollary: simply laced Gelfand-MacPherson resolution} may fail to hold.
Let \(W\) be the Weyl group of type \(C_2\) with Dynkin diagram
\[
\dynkin[label]{C}{2}
\]
and let \(w=s_2s_1s_2\).
It is well known that \(X_w\) is smooth.
This can be seen by taking the BP decomposition \(u_0=s_2s_1\) and \(u_1=s_2\) with respect to \(I=\tau(w)=\set{2}=\sigma(w)-\set{t}\), where \(t=s_1\).
Then \(X_{u_0}^{\tau(w)}\) is smooth by \cite[Theorem 3.8]{RichmondSlofstra16}, where we set \(W=C_2\), \(s=s_1\), and \(k=n=2\).
Hence \(X_w\) is also smooth, since it is a fiber bundle with base \(X_{u_0}^{\tau(w)}\) and fiber \(P_2/B\).
Observe that we have \(J=\sigma(u_0)\smallsetminus\set{t}=I\) and \(u_0\star w_J=w\), but \(w_{\sigma(u_0)}=s_2s_1s_2s_1\). 
It is clear that \(X_w\) does not admit resolution data such that the corresponding Gelfand-MacPherson resolution is an isomorphism.

Indeed, if there exists such an isomorphism, we can assume \(m=\ell-1\) (possibly with \(I_i=I_{i+1}\) for some \(i\)) such that for every \(0\leq i\leq m\), \(I_i\neq\emptyset\).
Let \(P_{I_0}\times^{R_1}P_{I_1}\times^{R_2}P_{I_2}\isom G_w\).
Note \(\#I_i=1\) since \(w<w_{\set{1,2}}\).
It follows that \(I_0=\set{2}=I_2\) and \(I_1=\set{1}\), which does not provide an isomorphism.
\end{example}

\begin{corollary}\label{simply laced tau}
Let \(X_w\) be a smooth simply laced Schubert variety.
Then \(\tau(w^{-1})=\tau(w)\) if and only if \(\tau(w)=\sigma(w)\).
\begin{proof}
If \(\tau(w)=\sigma(w)\) then \(w=w_{\tau(w)}\) by Lemma~\ref{lemma: sigma is tau}.
Hence \(w=w^{-1}\) in this case.

If \(\tau(w^{-1})=\tau(w)\), let \(\mu\,\colon P_{I_0}\times^{R_1}\cdots\times^{R_m}P_{I_m}/B\to X_w\) be an isomorphism such that \(I_0=\tau(w^{-1})\) and \(I_m=\tau(w)\) by Theorem~\ref{theorem: simply laced smooth isom}.
It follows that for every \(s\in\tau(w)\) and \(0\le i\le m\), we have \(s\in I_i\).
Indeed, if there exists \(0<i<m\) such that \(s\notin I_i\) then \([\dot s,1,\ldots,1,\dot sB/B]\) and \([1,\ldots,B/B]\) are different points in the fiber of \(\mu\) over \(B/B\).
This contradicts \(\mu\) being an isomorphism.
Hence for every \(0\le i\le m\), we have \(\tau(w)\subseteq I_i\).
Then for every \(0\le i\le m\), \(\tau(w)=I_i\) since \(I_m\) is always contained in \(\tau(w)\) (so \(I_m=\tau(w)\) in this case) and \(P_{I_0}\times^{P_{I_1}}P_{I_1}\isom P_{I_0}\) whenever \(I_1\subseteq I_0\).
Therefore \(\tau(w)=\sigma(w)\).
\end{proof}
\end{corollary}

\section{$A_{n-1}$}\label{section: A}

Fix \(G=GL(n,\bC)\) and let \(B\) be the upper triangular matrices in \(G\).
We recall a family of small resolutions described by Zelevinski\u\i\ \cite{Zelevinskii83}, and we use Lemma~\ref{lemma: small glue} to provide a new family of small resolutions in Proposition~\ref{proposition: small BP family}.
This family of small resolutions can be summarized using pattern avoidance.
Then we describe all Schubert varieties with small resolutions for \(A_{n-1}\) (\(n\leq 6\)).
We conclude with an example to show that pattern avoidance does not characterize the property `\(X_w\) admits a small resolution'.

\(G\) is of type \(A_{n-1}\) acting on the left of \(\bC^n\) as usual.
The standard basis \(\set{e_1,\ldots,e_n}\) of \(\bC^n\) fixes our choice of maximal torus \(T\subseteq B\) as the stabilizer of all lines \(\left\langle e_i\right\rangle\).
We identify the Weyl group \(W=N_G(T)/T\) with \(S_n\), the set of permutations of \(\set{1,\ldots,n}\), by letting \(\left\langle e_{w(i)}\right\rangle=\dot w\left\langle e_i\right\rangle\).
We denote a permutation \(w\) in one-line notation \(w=(\ w(1)\ \cdots\ w(n)\ )\).
The simple roots in the Dynkin diagram are labeled by
\[\dynkin[labels={1,2,,n-1}]{A}{}.\]

\begin{remark}
All resolutions in this section are Gelfand-MacPherson resolutions (as in \cite{GelfandMacPherson82} and \S\ref{section: Gelfand-MacPherson}).
The reason for this is explained in Corollary~\ref{corollary: simply laced Gelfand-MacPherson resolution}.
As a result, the resolutions can be described explicitly as an iterated base change, and a formula for fiber dimensions is provided by \cite{SankaranVanchinathan94}.
\end{remark}

Zelevinski\u\i\ \cite{Zelevinskii83} described a family of resolutions for every grassmannian Schubert variety for \(G\) by using a general construction of Gelfand-MacPherson (as described in \S\ref{section: Gelfand-MacPherson}).
He also showed each grassmannian Schubert variety has at least one small resolution.

Zelevinski\u\i\ used the iterated base change provided by Gelfand-MacPherson \cite{GelfandMacPherson82} to describe the resolutions in terms of incidence relations of flags.
Here we return to the description of resolutions using Bott-Samelson type varieties, following, e.g., \cite{SankaranVanchinathan94}, and the original construction of Demazure.

Let \(1\le k\le n\), \(\hat k=\set{1,\ldots,n}\smallsetminus\set{k}\), and consider a grassmannian Schubert variety \(X_w^{\hat k}\subseteq X^{\hat k}=G/P_{\hat k}\), where we are choosing \(w\) to be maximal in its \(W_{\hat k}\)-coset.
We point out that \(w\) maximal in its coset is equivalent to
\begin{equation*}\label{equation: grass perm}
w(1)>w(2)>\cdots>w(k),\quad w(k+1)>\cdots>w(n).
\end{equation*}

All of Zelevinski\u\i's resolutions (as mentioned) are
\begin{equation}\label{equation: Zelevinskii}
P_{I_0}\times^{R_1}\cdots\times^{R_m}P_{\hat k}/P_{\hat k}\to X_w^{\hat k}.
\end{equation}
It is important for us that in each of the resolutions of \cite{Zelevinskii83}, \(I_0=\tau(w^{-1})\). 
In the language of \cite{Zelevinskii83}, \(I_0=\set{s_j\mid j\text{ is not a valley}}\).
The \emph{valleys} are the \(j\neq n\) that begin each string of consecutive terms in \((w(1),\ldots,w(k))\).
As \(\tau(w^{-1})=\set{s_j\mid j+1\text{ appears left of }j\text{ in }w}\), we have that \(\tau(w^{-1})=S\smallsetminus\set{\text{valleys}}=I_0\).
For example, in type \(A_7\) with \(k=4\) and
\(w=\text{( 8 5 3 2 7 6 4 1 )}\),
the valleys are \(5\) and \(3\), and \(\tau(w^{-1})=\set{1,2,4,6,7}\).
When the resolutions \eqref{equation: Zelevinskii} are pulled back to resolutions of \(X_w\subseteq G/B\), they become
\begin{equation}\label{equation: Zelevinskii pull-back}
P_{I_0}\times^{R_1}\cdots\times^{R_m}P_{\hat k}/B\to X_w.
\end{equation}
When \(w\) is not equal to the long element of \(W\), \(\tau(w)=\hat k\).

This discussion shows that a restatement of the main result of \cite{Zelevinskii83} is the following.

\begin{theorem}[\cite{Zelevinskii83}]\label{theorem: Zelevinskii}
If \(w\in S_n\) is maximal in its \(W_{\hat k}\)-coset, then there is a small resolution
\begin{equation*}
P_{I_0}\times^{R_1}\cdots\times^{R_m}P_{I_m}/B\to X_w
\end{equation*}
satisfying \eqref{equation: equivariance}, i.e., \(I_0=\tau(w^{-1})\) and \(I_m=\tau(w)\).
\end{theorem}

\begin{corollary}\label{corollary: Zelevinskii}
If \(w\in S_n\) satisfies \(\#\tau(w)\geq \#\sigma(w)-1\), then there exists a small resolution \(\mu\,\colon P_{I_0}\times^{R_1}\cdots\times^{R_m}P_{I_m}/B\to X_w\) satisfying \eqref{equation: equivariance}.

\begin{proof}
Let \(w\in W=S_n\) such that \(\#\tau(w)\geq \#\sigma(w)-1\).
If \(\#\tau(w)=\#\sigma(w)\) then \(X_w=P_{\sigma(w)}/B\) by Lemma~\ref{lemma: sigma is tau}, and we are done.
We can assume that \(\sigma(w)\) is connected by applying Lemma~\ref{lemma: small glue} to \(G_{w_0}\times^B\cdots\times^B X_{w_m}\to X_w\), where \(\sigma(w_i)\) are pairwise disjoint and non-adjacent, so an isomorphism by Corollary~\ref{corollary: fiber bundle decomp}.
Then, for example, repeatedly applying \eqref{equation: parabolic inflation} gives the desired resolution satisfying \eqref{equation: equivariance}.

If \(\#\tau(w)=\#\sigma(w)-1\) then \(X_w^{\tau(w)}\) is isomorphic to a grassmannian Schubert variety for a smaller group of type \(A_{\#\sigma(w)}\) by Facts~\ref{facts: grassmannian} (b).
There exists resolution data for the corresponding grassmannian Schubert variety given by Corollary~\ref{corollary: Zelevinskii}.
The corresponding parabolic subgroups of the original \(G\) gives resolution data for \(X_w^{\tau(w)}\) by Proposition~\ref{proposition: mu}, since birational holds true by \eqref{equation: decomposition}. 
The corresponding resolution is small since the formula for fiber dimensions in \cite{SankaranVanchinathan94} shows the dimensions are the same.
We have \(I_0=\tau(w^{-1})\) since this holds true for the resolution in the smaller group.
By Theorem~\ref{theorem: Zelevinskii}, we have a small resolution of \(X_w\) with \(I_0=\tau(w^{-1})\) and \(I_m=\tau(w)\).
\end{proof}
\end{corollary}

\begin{example}
Let \(w=\text{( 4 2 3 1 )}\) with reduced expression \(w=s_1s_3s_2s_1s_3\).
In this case \(\tau(w^{-1})=\set{1,3}=\tau(w)\).
Then \(\#\tau(w)=\#\sigma(w)-1\) (and \(w\) satisfies the hypothesis of Theorem~\ref{theorem: Zelevinskii}), so \(X_w\) has a small resolution by Corollary~\ref{corollary: Zelevinskii} (and Theorem~\ref{theorem: Zelevinskii}).
By Theorem~\ref{theorem: Zelevinskii}, the two small resolutions corresponding to `neat ordering of peaks', as defined in \cite{Zelevinskii83}, can be described by \(\mu\,\colon P_{1,3}\times^{P_3}P_{2,3}\times^{P_3}P_{1,3}/B\to X_w\) and \(\nu\,\colon P_{1,3}\times^{P_1}P_{1,2}\times^{P_1}P_{1,3}/B\to X_w\).
\end{example}

\begin{example}
Let \(w=\text{( 1 5 3 4 2 )}\).
Note \(X_w\) is not the pull-back of a grassmannian Schubert variety, but \(\#\tau(w)=\#\sigma(w)-1\), so is isomorphic to the pull-back of a grassmannian Schubert variety \(X_u\) for a smaller group, where \(u=\text{( 4 2 3 1 )}\).
\end{example}

In this section, we obtain a new family of small resolutions by applying Lemma~\ref{lemma: small glue} to \cite{Zelevinskii83}.
The family is best described by recalling a pattern avoidance result of \cite{AllandRichmond18}.
Then using Proposition~\ref{proposition: inverse trick}, we see that the family extends to be stable under the function \(w\mapsto w^{-1}\).

\begin{proposition}\label{proposition: small BP family}
If \(w\) avoids the patterns
\begin{equation*}\label{equation: complete BP patterns}
\text{( 3 4 1 2 ), ( 5 2 3 4 1 ), ( 6 3 5 2 4 1 )}
\end{equation*}
then \(X_w\) and \(X_{w^{-1}}\) have small resolutions.

\begin{proof}
\cite[Theorem 1.4, Proposition 2.6]{AllandRichmond18} shows that \(w\) avoids this list of patterns if and only if it has a complete BP decomposition.
Hence we can apply Facts~\ref{facts: complete BP decomp} to get a fiber bundle decomposition \(G_{w_0}\times^{R_1}\cdots\times^{R_m}X_{w_m}\to X_w\) such that for every \(0\le i\le m\), we have \(\tau(w_i)=\sigma(w_i)\) or \(\tau(w_i)=\sigma(w_i)\smallsetminus\set{s_i}\).
By Corollary~\ref{corollary: Zelevinskii}, for every \(0\le i\le m\), \(X_{w_i}\) admits a small resolution satisfying \eqref{equation: equivariance}.
Hence we can use Lemma~\ref{lemma: small glue} to obtain a small resolution of \(X_w\).
Then Proposition~\ref{proposition: inverse trick} gives us a small resolution of \(X_{w^{-1}}\).
\end{proof}
\end{proposition}

\begin{example}
Let \(w=\text{( 6 3 5 2 4 1 )}\).
Note \(w\) does not satisfy Proposition~\ref{proposition: small BP family}.
Then \(w^{-1}=\text{( 6 4 2 5 3 1 )}\) satisfies Proposition~\ref{proposition: small BP family} (and Corollary~\ref{corollary: Zelevinskii}).
Therefore \(X_w\) has a small resolution.
\end{example}

\begin{example}
Let \(w=\text{( 6 4 5 7 3 2 1 )}\).
Then \(w\) satisfies Proposition~\ref{proposition: small BP family}.
The decomposition \[\tilde w=(s_3 s_2 s_1 s_5 s_4 s_3 s_2,s_1,s_5 s_4 s_3,s_6 s_5 s_4,s_6 s_5,s_6)\] is a complete BP decomposition.
As in \eqref{equation: w_i}, let
\begin{equation*}
\begin{split}
w_0&=s_1 s_2 s_1 s_3 s_2 s_1 s_4 s_3 s_2 s_5 s_4 s_3 s_2 s_1 =\text{( 6 4 5 3 2 1 7 )}\\ 
w_1&=s_1 s_3 s_4 s_3 s_5 s_4 s_3 =\text{( 2 1 6 5 4 3 7 )}\\ 
w_2&=s_3 s_4 s_3 s_5 s_4 s_3 s_6 s_5 s_4 =\text{( 1 2 6 7 5 4 3 )}\\ 
w_3&=s_4 s_5 s_4 s_6 s_5 s_4 =\text{( 1 2 3 7 6 5 4 )}\\ 
w_4&=s_5 s_6 s_5 =\text{( 1 2 3 4 7 6 5 )}\\ 
w_5&=s_6 =\text{( 1 2 3 4 5 7 6 )}\\ 
\end{split}
\end{equation*}
with corresponding isomorphism \(\mu\,\colon G_{w_0}\times^{R_1}\cdots\times^{R_5}X_{w_5}\to X_w\).

Then 
\(P_{\set{1,2,3,5}}\times^{P_{\set{1,2,3}}}P_{\set{1,2,3,4}}\times^{P_{\set{1,3,4}}}P_{\set{1,3,4,5}}\to G_{w_0}\)
is a small resolution such that \(\tau(w_0^{-1})=I_0\) and \(\tau(w_0)=I_2\).
For \(1\le i\le 5\), \(X_{w_i}\) is smooth since \(w_i\) avoids \(3412\) and \(4231\).
Hence \(\mu\,\colon P_{\set{1,2,3,5}}\times^{P_{\set{1,2,3}}}P_{\set{1,2,3,4}}\times^{P_{\set{1,3,4}}}P_{\set{1,3,4,5}}\times^{R_1}G_{w_1}\times^{R_2}\cdots\times^{R_5}X_{w_5}\to X_w\) is a small resolution.
\end{example}

We provide an example in \S\ref{example: no pattern avoidance} to show that the property `\(X_w\) admits a small resolution' is not characterized by pattern avoidance.
Along the way we provide data to show which Schubert varieties admit small resolutions in \(W=S_5\) of type \(A_4\) and \(W=S_6\) of type \(A_5\).
We conclude that for \(n\le 6\) and \(w\in W=S_{n}\) of type \(A_{n-1}\), then \(X_w\) has a small resolution if and only if \(X_w\) does not have factorial singular locus.

Let \(W=S_5\) of type \(A_4\).
There are 120 Schubert varieties in \(X\), and 119 of these have small resolutions.
The remaining Schubert variety corresponding to \(w=\text{( 4 5 3 1 2 )}\) is known to be singular and factorial by \cite{WooYong08}.
It is well-known that a singular and factorial (or more generally \(\bQ\)-factorial) algebraic variety does not admit \emph{any} small resolution (as e.g., in \cite{Perrin07}).

There are 88 smooth Schubert varieties, so the small resolutions in this case are the identity morphism.
There are 8 singular Schubert varieties with small resolutions by \cite{BilleyWarrington01} (avoiding 321-hexagon patterns) and 14 by Proposition~\ref{proposition: small BP family}.
Table~\ref{table: A4} provides a description for small resolutions of the form \(P_{\tau(w^{-1})}\times^{R_1}G_{w_1}\times^{R_2}P_{\tau(w)}\to X_w\) for the remaining 9 singular Schubert varieties with small resolutions.
This table was constructed by finding \(w_1\) such that \(w=w_{I_0}\star w_1\star w_{I_2}\), where \(I_0=\tau(w^{-1})\), \(I_2=\tau(w)\), and the dimension formula of \cite{SankaranVanchinathan94} shows smallness.
This was accomplished with help of the atlas software \cite{atlas}.

\begin{table}[h!]
\caption{Small resolutions for \(W=S_5\)}
\begin{tabular}{c|rcl|ll}
$w$ & \(\tau(w^{-1})\) & \(w_1\) & \(\tau(w)\) & \(\tau(w_1^{-1})\) & \(\tau(w_1)\)\\
\hline
( 3 5 1 4 2 ) & $\set{2,4}$   & ( 2 1 5 4 3 ) & $\set{2,4}$   & $\set{1,3,4}$ & $\set{1,3,4}$ \\
( 4 2 5 1 3 ) & $\set{1,3}$   & ( 3 2 1 5 4 ) & $\set{1,3}$   & $\set{1,2,4}$ & $\set{1,2,4}$ \\
( 4 5 1 3 2 ) & $\set{2,3}$   & ( 2 1 5 4 3 ) & $\set{2,4}$   & $\set{1,3,4}$ & $\set{1,3,4}$ \\
( 3 5 4 1 2 ) & $\set{2,4}$   & ( 2 1 5 4 3 ) & $\set{2,3}$   & $\set{1,3,4}$ & $\set{1,3,4}$ \\
( 4 3 5 1 2 ) & $\set{2,3}$   & ( 3 2 1 5 4 ) & $\set{1,3}$   & $\set{1,2,4}$ & $\set{1,2,4}$ \\
( 4 5 2 1 3 ) & $\set{1,3}$   & ( 3 2 1 5 4 ) & $\set{2,3}$   & $\set{1,2,4}$ & $\set{1,2,4}$ \\
( 5 2 3 4 1 ) & $\set{1,4}$   & ( 1 4 3 2 5 ) & $\set{1,4}$   & $\set{2,3}$   & $\set{2,3}$ \\
( 5 3 4 1 2 ) & $\set{2,4}$   & ( 4 3 1 5 2 ) & $\set{1,3}$   & $\set{2,3}$   & $\set{1,2,4}$ \\
( 4 5 2 3 1 ) & $\set{1,3}$   & ( 4 1 5 3 2 ) & $\set{2,4}$   & $\set{2,3}$   & $\set{1,3,4}$ \\
\end{tabular}
\label{table: A4}
\end{table}

A similar classification holds for \(W=S_6\) of type \(A_5\).
There are 720 Schubert varieties in \(X\), and exactly 701 of these have small resolutions.

There are 366 smooth Schubert varieties, 43 singular Schubert varieties \(X_w\) such that \(w\) avoids 321-hexagon patterns, and 127 singular Schubert varieties satisfying Proposition~\ref{proposition: small BP family} (55 for which \(w\) or \(w^{-1}\) satisfy Corollary~\ref{corollary: Zelevinskii}).
Out of the remaining 165 Schubert varieties with small resolutions, 56 have fiber bundle decompositions \(X_u\isom G_{v}\times^{R}X_{w}\) such that \(v,w<u\).
We remark that Proposition~\ref{proposition: small BP family} does not assert that the small resolution satisfies \eqref{equation: equivariance} (so care must be taken when applying Lemma~\ref{lemma: small glue}), but we have checked that this does hold true for \(n\le6\).

There are 109 Schubert varieties with small resolutions that are not described by above considerations, and 91 of these \(X_w\) have the property that \(\#\sigma(w)=5\).
These resolutions were found using atlas software \cite{atlas} to compute fiber dimensions of Gelfand-MacPherson resolutions.
One can find many small resolutions recursively by first looking for small resolutions satisfying \eqref{equation: equivariance}.
We provide in Table ~\ref{table: A5}, 53 Schubert varieties \(X_w\) such that all \(w\) or \(w^{-1}\) provides the list of 91 small resolutions above.
To reconstruct the small resolution from Table ~\ref{table: A5}, let \((I_0,\ldots,I_m)\) give a small resolution of \(X_{w_1}\) such that \(I_0=\tau(w_1^{-1})\) and \(I_m=\tau(w_1)\).
Then \((\tau(w^{-1}),I_0,\ldots,I_m,\tau(w))\) gives a small resolution of \(X_w\).
This accounts for all Schubert varieties having small resolutions.

There are 19 Schubert varieties that are either singular and factorial, or contain the (singular and factorial) interval \([\text{ 1 4 3 2 5 },\text{ 4 5 3 1 2 }]\).
It follows that these Schubert varieties do not admit \emph{any} small resolution.

\begin{example}\label{example: no pattern avoidance}
Let \(w=\text{( 4 6 3 1 5 2 )}\) in \(W=S_6\) of type \(A_5\), so \(\tau(w^{-1})=\set{2,3,5}=\tau(w)\).
Let \(I_0=\set{2,3,5}\), \(I_1=\set{1,2,4,5}\), and \(I_2=\set{2,3,5}\).
Then \(\mu\,\colon P_{I_0}\times^{R_1}P_{I_1}\times^{R_2}P_{I_2}/B\to X_w\) is a small resolution by Table ~\ref{table: A5}, where \(J_1=\set{2,5}=J_2\).
The permutation \(w\) contains the pattern \(u=45312\), and \(X_u\) does not have a small resolution since it is factorial.
Therefore small resolutions are not characterized by pattern avoidance.
\end{example}

\begin{table}[h!]
\caption{Small resolutions for $W=S_6$}
\begin{minipage}{.4\textwidth}
\begin{tabular}{c|c}
$w$ & $w_1$\\
\hline
( 4 6 1 2 5 3 ) & ( 3 1 6 2 5 4 ) \\
( 3 6 1 4 5 2 ) & ( 2 1 5 4 3 6 ) \\
( 5 2 6 1 3 4 ) & ( 4 2 1 6 3 5 ) \\
( 4 2 6 1 5 3 ) & ( 3 2 1 6 5 4 ) \\
( 5 2 3 6 1 4 ) & ( 1 4 3 2 6 5 ) \\
( 5 6 1 2 4 3 ) & ( 3 1 6 2 5 4 ) \\
( 4 6 1 5 2 3 ) & ( 3 1 6 5 2 4 ) \\
( 5 6 1 3 2 4 ) & ( 4 1 6 3 2 5 ) \\
( 4 6 1 3 5 2 ) & ( 2 1 5 4 3 6 ) \\
( 5 3 6 1 2 4 ) & ( 4 3 1 6 2 5 ) \\
( 3 6 1 5 4 2 ) & ( 2 1 6 5 4 3 ) \\
( 4 3 6 1 5 2 ) & ( 3 2 1 6 5 4 ) \\
( 4 3 5 6 1 2 ) & ( 3 2 5 1 6 4 ) \\
( 5 2 6 1 4 3 ) & ( 3 2 1 6 5 4 ) \\
( 5 2 4 6 1 3 ) & ( 1 4 3 2 6 5 ) \\
( 6 2 3 4 5 1 ) & ( 1 5 3 4 2 6 ) \\
( 5 3 2 6 1 4 ) & ( 4 3 2 1 6 5 ) \\
( 4 6 5 1 2 3 ) & ( 3 1 6 5 2 4 ) \\
( 5 4 6 1 2 3 ) & ( 4 3 1 6 2 5 ) \\
( 5 6 1 3 4 2 ) & ( 4 1 6 3 2 5 ) \\
( 4 6 1 5 3 2 ) & ( 2 1 6 5 4 3 ) \\
( 5 3 6 1 4 2 ) & ( 4 3 1 6 5 2 ) \\
( 6 3 4 1 5 2 ) & ( 5 2 1 4 3 6 ) \\
( 5 3 4 6 1 2 ) & ( 2 5 4 1 6 3 ) \\
( 4 6 3 1 5 2 ) & ( 3 2 1 6 5 4 ) \\
( 4 3 6 5 1 2 ) & ( 3 2 1 6 5 4 ) \\
( 6 2 4 5 1 3 ) & ( 1 5 4 2 6 3 ) \\
\end{tabular}
\end{minipage}
\begin{minipage}{.4\textwidth}
\begin{tabular}{c|c}
$w$ & $w_1$\\
\hline
( 5 2 6 4 1 3 ) & ( 3 2 1 6 5 4 ) \\
( 5 4 2 6 1 3 ) & ( 4 3 2 1 6 5 ) \\
( 6 2 3 5 4 1 ) & ( 1 6 3 5 4 2 ) \\
( 6 2 4 3 5 1 ) & ( 1 5 4 3 2 6 ) \\
( 6 3 2 4 5 1 ) & ( 5 3 2 4 1 6 ) \\
( 6 4 5 1 2 3 ) & ( 5 4 1 6 2 3 ) \\
( 5 6 1 4 3 2 ) & ( 2 1 6 5 4 3 ) \\
( 4 6 5 1 3 2 ) & ( 2 1 6 5 4 3 ) \\
( 5 4 6 1 3 2 ) & ( 3 2 1 6 5 4 ) \\
( 6 3 5 1 4 2 ) & ( 6 2 1 5 4 3 ) \\
( 5 6 3 1 4 2 ) & ( 4 3 1 6 5 2 ) \\
( 6 3 4 5 1 2 ) & ( 2 5 4 1 6 3 ) \\
( 5 3 6 4 1 2 ) & ( 4 3 1 6 5 2 ) \\
( 5 4 3 6 1 2 ) & ( 4 3 2 1 6 5 ) \\
( 6 4 2 5 1 3 ) & ( 5 4 2 1 6 3 ) \\
( 5 4 6 2 1 3 ) & ( 4 3 2 1 6 5 ) \\
( 6 2 5 3 4 1 ) & ( 1 5 4 3 2 6 ) \\
( 6 4 2 3 5 1 ) & ( 1 5 4 3 2 6 ) \\
( 6 4 5 1 3 2 ) & ( 6 2 1 5 4 3 ) \\
( 6 3 5 4 1 2 ) & ( 2 6 5 4 1 3 ) \\
( 5 6 3 4 1 2 ) & ( 4 6 3 1 5 2 ) \\
( 6 4 3 5 1 2 ) & ( 5 4 3 1 6 2 ) \\
( 6 4 5 2 1 3 ) & ( 5 4 2 1 6 3 ) \\
( 6 5 2 3 4 1 ) & ( 1 5 4 3 2 6 ) \\
( 6 5 3 4 1 2 ) & ( 5 4 3 1 6 2 ) \\
( 6 4 5 2 3 1 ) & ( 5 4 1 6 3 2 ) \\
&\\
\end{tabular}
\end{minipage}

\label{table: A5}
\end{table}





\newpage
\bibliography{../../references/my}{}

\end{document}